\documentclass[journal]{IEEEtran}
\usepackage{amsmath}
\usepackage[ruled]{algorithm2e}
\usepackage{amssymb,amsthm,bm,url,graphicx}

\usepackage{jdhmacros}
\usepackage{setspace}
\usepackage{bbm}
\usepackage{dsfont}
\usepackage{multirow}
\usepackage{xcolor}
\usepackage{tabularx, booktabs}
\usepackage{multirow}
\usepackage{epsfig}
\usepackage{cite}

\newcommand{\st}{\text{s.t.}}

\newcommand{\jdh}[1]{}
\renewcommand{\jdh}[1]{{\color{magenta}{#1}}} 
\ifCLASSINFOpdf
\else
\fi
  \usepackage[caption=false,font=footnotesize]{subfig}
\begin{document}
%
\righthyphenmin = 2
\lefthyphenmin = 2
\title{A Provably Communication-Efficient\\ Asynchronous Distributed Inference Method\\ for Convex and Nonconvex Problems}
%
%
%

\author{Jineng Ren,~\IEEEmembership{Student Member,~IEEE,} and 
        Jarvis Haupt,~\IEEEmembership{Senior~Member,~IEEE}%
\thanks{The authors are with the Department
of Electrical and Computer Engineering, University of Minnesota, Minneapolis,
MN, 55455 USA; e-mails: \{renxx282,jdhaupt\}@umn.edu.  Shorter preliminary versions of this work appeared at the 2018 Global Conference on Signal and Information Processing (GlobalSIP 2018).}
\thanks{\textcopyright20XX IEEE.  Personal use of this material is permitted.  Permission from IEEE must be obtained for all other uses.}}
\maketitle

\begin{abstract} 
This paper proposes and analyzes a communication-efficient distributed optimization framework for general nonconvex nonsmooth signal processing and machine learning problems under an asynchronous protocol. At each iteration, worker machines compute gradients of a known empirical loss function using their own local data, and a master machine solves a related minimization problem to update the current estimate. We prove that for nonconvex nonsmooth problems, the proposed algorithm converges with a sublinear rate over the number of communication rounds, coinciding with the best theoretical rate that can be achieved for this class of problems. Linear convergence is established without any statistical assumptions of the local data for problems characterized by composite loss functions whose smooth parts are strongly convex. Extensive numerical experiments verify that the performance of the proposed approach indeed improves -- sometimes significantly -- over other state-of-the-art algorithms in terms of total communication efficiency.
\end{abstract}


%
\IEEEpeerreviewmaketitle

\section{Introduction}
%
%
%
%
\IEEEPARstart{D}{ue} to rapid developments in information and computing technology, modern applications often involve vast amounts of data, rendering local processing (e.g., in a single machine, or on a single processing core) computationally challenging or even prohibitive. To deal with this problem, distributed and parallel implementations are natural methods that can fully leverage multi-core computing and storage technologies. However, one drawback of distributed algorithms is that the communication cost can be very expensive in terms of raw bytes transmitted, latency, or both, as machines (i.e., computation nodes) need to frequently transmit and receive information between each other. Therefore, algorithms that require less communication are preferred in this case. 

In this paper we study a general communication-efficient distributed algorithm which can be applied to a broad class of nonconvex nonsmooth inference problems. Assume that we have available some $N$ data samples. 
We consider a general problem appearing frequently in signal processing and machine learning applications; we aim to solve
\begin{align}\label{origin0}
\xb^{*} = \argmin_{\xb\in \RR^p} \Lb(\xb) :=\frac{1}{N}\sum_{k=1}^N l_{k}(\xb) + h(\xb), 
\end{align}
where each $l_{k}(\xb)$ is a loss function associated with the $k$-th data sample, and is assumed smooth but possibly nonconvex with Lipschitz continuous gradient and $h(\xb)$ is a convex (proper and lower semi-continuous) function that is possibly nonsmooth. Problem \eqref{origin} covers many important machine learning and signal processing problems such as the localization with wireless acoustic sensor networks (WASNs) \cite{meguerdichian2001localized}, support vector machine (SVM) \cite{friedman2001elements}, the independent principal component analysis (ICA) reconstruction problem \cite{le2011ica}, and the sparse principal component analysis (PCA) problem \cite{richtarik2012alternating}.

\begin{figure}[!tbh]
	\centering 
	\vspace{0.01in}
	\begin{center}
		\hspace{-0.096in}{
			\includegraphics[width=0.9\linewidth]{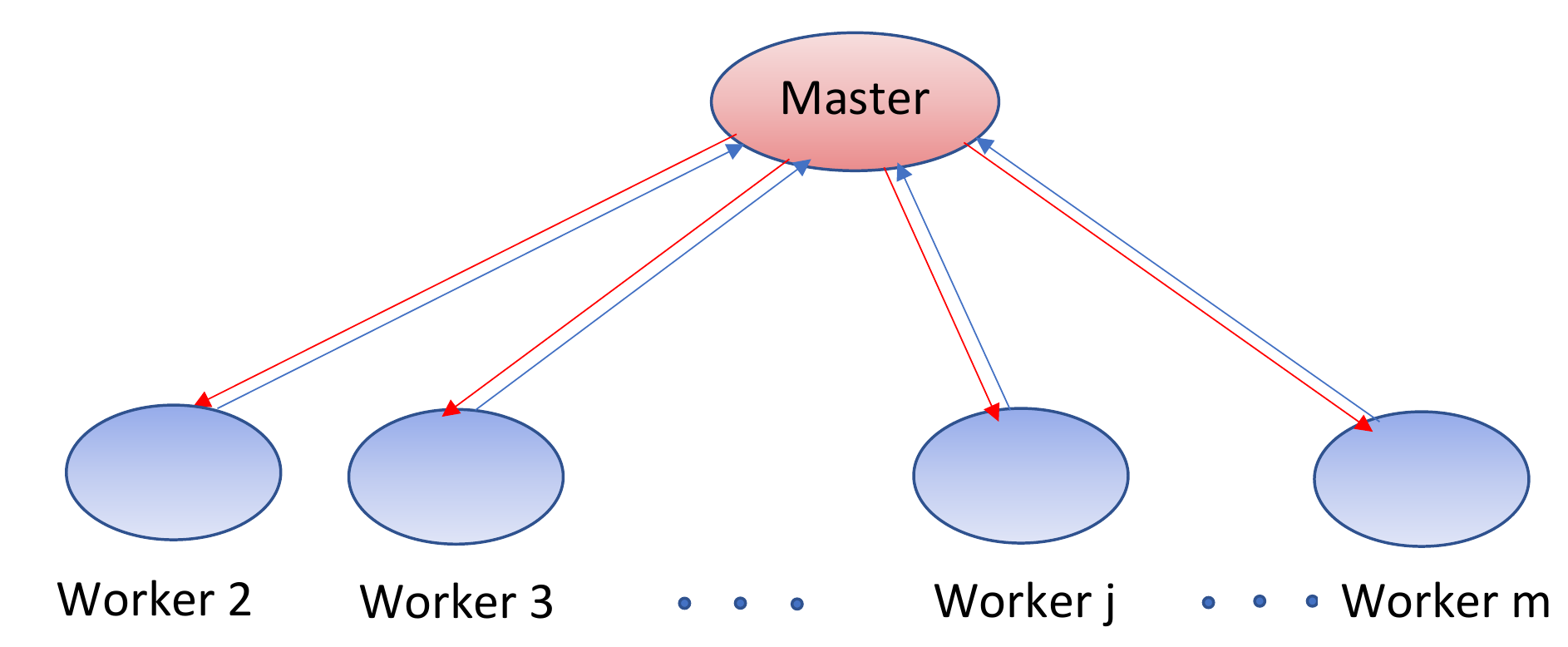}
		}
	\end{center}
	\vspace{-0.1in}	
	\caption{$m$-nodes network with a star topology}\label{networkfig}
	\vspace{-0.12in}
\end{figure}

For our distributed approach, we consider a network of $m$ total machines having a star topology, where one node designated as the ``Master'' node (node $1$, without loss of generality) is located at the center of the star, and the remaining $m-1$ nodes (with indices $2,3,\dots,m$) are the ``Worker'' nodes (see Figure~\ref{networkfig}). Without loss of generality, assume that the number of data samples is evenly divisible by $m$, i.e., $N = nm$ for some integer $n$, and each machine stores $n$ unique data samples. Then \eqref{origin0} can be reformulated to the following problem: 
\begin{align}\label{origin}
\xb^{*} = \argmin_{\xb\in \RR^p} \Lb(\xb) :=\frac{1}{mn}\sum_{j=1}^m \sum_{i =1}^n l_{ji}(\xb) + h(\xb), 
\end{align}
where $l_{ji}(\xb)$ is the loss function corresponding to the $i$-th sample of the $j$-th machine. 

\subsection{Main Results}
We propose an Efficient Distributed Algorithm for Nonconvex Nonsmooth Inference (EDANNI), and show that, for general problems of the form of \eqref{origin}, EDANNI converges to the set of stationary points if the algorithm parameters are chosen appropriately according to the maximum network delay. Our results differ significantly from existing works \cite{rabbat2004distributed,agarwal2011distributed,bradley2011parallel} which are all developed for convex problems. Therefore, the analysis and algorithm proposed here are applicable not only to standard convex learning problems but also to important nonconvex problems. 
To the best of our knowledge, this is the first communication-efficient algorithm exploiting local second-order information that is guaranteed to be convergent for general nonconvex nonsmooth problems. Moreover, linear convergence is also proved in the strongly convex setting with no statistical assumption of the data stored in each local machine, which is another improvement on existing works. The synchronization inherent in previous works, including \cite{rabbat2004distributed,agarwal2011distributed,bradley2011parallel,pmlr-v70-wang17f,jordan2016communication,ren2017communication}, slows down those methods because the master needs to wait for the slowest worker during each iteration; here, we propose an asynchronous approach that can accelerate the corresponding inference tasks significantly (as we will demonstrate in the experimental results). 

\subsection{Related Work}  
There is a large body of work on distributed optimization for modern data-intensive applications with varied accessibility; see, for example, \cite{rabbat2004distributed,agarwal2011distributed,bradley2011parallel,tsitsiklis1986distributed,recht2011hogwild,hong2014distributed,boyd2011distributed,zinkevich2010parallelized,bertsekas1989parallel,dean2012large,li2014scaling,li2014communication,chang2016asynchronous,chang2016asynchronous1,ren2018b,ren2018c}. (Parts of the results presented here appeared in our conference paper \cite{ren2018b} without theoretical analysis.)  Early works including \cite{agarwal2011distributed,tsitsiklis1986distributed,zinkevich2010parallelized} mainly considered the convergence of parallelizing stochastic gradient descent schemes which stem from the idea of the seminal text by Bertsekas and Tsitsiklis \cite{bertsekas1989parallel}. Niu et al. \cite{recht2011hogwild} proposed a lock-free implementation of distributed SGD called Hogwild! and provided its rates of convergence for sparse learning problems. That was followed up by many variants like \cite{cannelli2016asynchronous, liu2015asynchronous}. For solving large scale problems, works including \cite{dean2012large}, \cite{li2014scaling,li2014communication}, and \cite{aytekin2016analysis} studied distributed optimizations based on a parameter server framework and parameters partition. Chang et al. \cite{chang2016asynchronous} studied asynchronous distributed optimizations based on the alternating direction method of multipliers (ADMM). By formulating the optimization problem as a consensus problem, the ADMM can be used to solve the consensus problem in a fully parallel fashion over networks with a star topology. One drawback of such approaches is that they can be computationally intensive, since each worker machine is required to solve a high dimensional subproblem. As we will show, these methods also converge more slowly (in terms of communication rounds) as compared to the proposed approach (see Section \ref{experimentsection}). 


A growing interest on distributed algorithms also appears in the statistics community \cite{dekel2012optimal,zhang2012communication,shamir2014distributed,arjevani2015communication,lee2015distributed}. Most of these algorithms depend on the partition of data, so their work usually involves statistical assumptions that handle the correlation between the data in local machines. A popular approach in early works is averaging estimators generated locally by different machines \cite{zinkevich2010parallelized,zhang2012communication,mcdonald2009efficient, huang2015distributed}. Yang \cite{yang2013trading},  Ma et al. \cite{ma2015adding}, and Jaggi et al. \cite{jaggi2014communication} studied distributed optimization based on stochastic dual coordinate descent, however, their communication complexity is not better than that of first-order approaches. Shamir et al. \cite{shamir2014communication} and Zhang and Xiao \cite{zhang2015disco} proposed truly communication-efficient distributed optimization algorithms which leveraged the local second-order information, though these approaches are only guaranteed to work for convex and smooth objectives. In a similar spirit, Wang et al. \cite{pmlr-v70-wang17f}, Jordan et al. \cite{jordan2016communication}, and Ren et al. \cite{ren2017communication} developed communication-efficient algorithms for sparse learning with $\ell_1$ regularization. However, each of these works needs an assumption about the strong convexity of loss functions, which may limit their approaches to only a small set of real-world applications. Here we describe an algorithm with similar flavor, but with more general applicability, and establish its convergence rate in both strongly convex and nonconvex nonsmooth settings. Moreover, unlike \cite{pmlr-v70-wang17f,jordan2016communication,ren2017communication,shamir2014communication,zhang2015disco} where the convergence analyses rely on certain statistical assumptions on the data stored in machines, our convergence analysis is deterministic and characterizes the worst-case convergence
conditions.

\hspace{-1.05em}\textbf{Notation.}
For a vector $\vb = (v_1,\cdots, v_s)^\top \in \mathbb{R}^s$ and $q>0$ we write $\|\vb\|_q = ( \sum_{i =1}^s |v_i|^q)^{1/q}$; for $q\geq1$ this is a norm. Usually $\|\vb\|_2$ is briefly written as $\|\vb\|$. The set of natural numbers is denoted by $\NN$. For an integer $m\in \NN$, we write $[m]$ as shorthand for the set $\{1,\dots,m\}$.

\section{Algorithm} 
In this section, we describe our approach to computing the minimizer $\xb^*$ of \eqref{origin}. Recall that we have $m$ machines. Let us denote $t\geq 0$ as the iteration number, then $\cA_t \subseteq [m] := \{1,2,\cdots,m\}$ is defined as the index of a subset of worker machines from which the master receives updated gradient information during iteration $t$; worker $i$ is said to be ``arrived" if $i\in \cA_t$. At iteration $t$, the master machine solves a subproblem to obtain an updated estimate, and communicates this to the worker machines in the subset $\cA_t$. After receiving the updated estimate, the worker machines will compute the corresponding gradients of local empirical loss functions. These gradients are then communicated back to the master machine, and the process continues. 

Formally, let 
\begin{align*}
\Lb_j(\xb) = \frac{1}{n}\sum_{i\in [n]}l_{ji}(\xb),~~j\in [m]
\end{align*}
be the empirical loss at each machine.  Let $t_j$ be the latest time (in terms of the iteration count) when the worker $j$ is arrived
up to and including iteration $t$. 

In the $t$-th iteration, the master (machine~$1$) solves the following subproblem to update $\xb^{t+1}$
\begin{align}\label{update}
\nonumber \xb^{t+1} &= \argmin_{\xb} \Lb_1(\xb) + h(\xb) + \frac{\rho}{2}\|\xb - \xb^t\|^2  \\ 
& + \Big\langle\frac{1}{m}\sum_{j \in [m]} \nabla \Lb_j(\xb^{t_j}) - \nabla\Lb_1(\xb^{t_1}), \xb - \xb^t \Big\rangle. 
\end{align}
This $\xb^{t+1}$ is communicated to the worker machines that are free, where it is used to compute their local gradient $\nabla\Lb_j(\xb^{t+1})$. Since machine $1$ is assumed to be the master machine, $t_1$ is actually $t$. 

Now one question is: which partial sets of worker machines (with indices in $\cA_t$) from which the master receives updated gradient information during iteration $t$ are sufficient to ensure convergence of a distributed approach? Firstly, let $\tau \geq 0$ be a maximum tolerable delay, that is, the maximum number of iterations for which every worker machine can be inactive. The set $\cA_t$ should satisfy:  

\begin{assumption}[Bounded delay]\label{boundeddelay}
	For all $i\in [m] $ and iteration $t\geq0$, it holds that $i\in \cA_t\cup\cA_{t-1}\cup\cdots\cup\cA_{\max\{t-\tau,0\}}$. 	
\end{assumption}

To satisfy Assumption \ref{boundeddelay}, $\cA_t$ should contain at least the indices of the worker machines that have been inactive for longer than $\tau$ iterations. That is, the master needs to wait until those workers finish their current computation and have arrived. Note that by the definition of $t_j$, it holds that 
\begin{align*}
t-\tau\leq t_j \leq t, \hspace{1.1em}  \forall j\in [m].    
\end{align*} 

Assumption~\ref{boundeddelay} requires that every worker $j$ is arrived at least once within the period $[t - \tau , t]$.
In other words, the gradient information $\nabla\Lb_j(\xb^{t_j})$ used by the master must be at most $\tau$ iterations old.
To guarantee the bounded delay, at every iteration the master needs to wait for the workers who have not been active for $\tau$ iterations, if such workers exist. Note that, when $\tau = 0$, one has $j\in \cA_t$ for all
$j\in[m]$, which reduces to the synchronous case where the master always waits for all the workers at every iteration.  

The proposed approach is presented in Algorithm~\ref{algor}, which specifies respectively the steps for the workers and the master. Algorithm~\ref{algor} has three prominent differences compared with its synchronous alternatives. First, only the workers $j$ in $\cA_t$ update the gradient $\nabla \Lb_j(\xb^{t})$ and transmit it to the master machine. For the workers $j$ in $\cA_t^c$, the master uses their latest gradient information before $t$, i.e.,  $\nabla \Lb_j(\xb^{t_j})$. Second, the variables $d_j$'s are introduced to count the delays of the workers since their last updates. $d_j$ is set to zero if worker $j$ is arrived at the current iteration; otherwise, $d_j$ is increased by one. Therefore, to ensure Assumption~\ref{boundeddelay} holds at each iteration, the master should wait if there exists at least one worker whose $d_j > \tau-1$. Third, after solving subproblem \eqref{update}, the master transmits the up-to-date variable $\xb^{t+1}$ only to the arrived workers. In general both the master and fast workers in the asynchronous approach can update more frequently and have less waiting time than their synchronous counterparts.

\begin{algorithm}[!htb]
	\caption{Efficient Distributed Algorithm for Nonconvex-Nonsmooth Inference (EDANNI)}
	\label{algor}
	\begin{algorithmic}
		\REQUIRE Loss functions $\{l_{ji}(\cdot,\cdot)\}_{i\in[n],j\in[m]}$, parameter $\rho$, 
		\STATE \hspace{2em} initial point $\xb^0$. Set $d_1= \cdots = d_m = 0$ and 
		\STATE \hspace{2em} $\cA_0 = [m]$; 
		\FOR{$t=0,1, \dots$}
		\STATE\textbf{\underline{Worker machines:}}
		\FOR{$j=2,3, \dots, m$}
		\STATE   \textbf{if} Receive $\xb^t$ from the master \textbf{then}
		\STATE Calculate gradient $\nabla \Lb_j(\xb^t)$ and transmit it to the master.
		\STATE \textbf{end}
		
		\ENDFOR
		\STATE  \textbf{\underline{Master:}}
		\STATE   \textbf{Receive} $\{\nabla\Lb_j(\xb^t)\}_{j=2}^m$ from worker machines $j$ in a set $\cA_t$ such that $d_j\leq\tau-1$,\hspace{0.5em} $\forall j \in \cA_t^c$
		\STATE \textbf{then}
		\STATE \hspace{1em}Update $$ d_j=\left\{
		\begin{aligned}
		&0 & \hspace{1.1em}\forall j\in \cA_t \\
		&d_j+1 & \hspace{1.1em}\forall j\in \cA_t^c
		\end{aligned} \hspace{0.5em}.
		\right.
		$$
		\STATE \hspace{1em}Solve the subproblem \eqref{update} with the specified $\rho$ 
		\STATE \hspace{1em}to obtain $\xb^{t+1}$. Broadcast $\xb^{t+1}$ to the worker 
		\STATE \hspace{1em}machines $j$ that are free. 
		\ENDFOR	
	\end{algorithmic}
\end{algorithm}

\section{Theoretical Analysis}
Solving subproblem $\eqref{update}$ is inspired by the approaches of Shamir et al. \cite{shamir2014communication}, et al., Wang et al. \cite{pmlr-v70-wang17f}, and Jordan et al.\cite{jordan2016communication}, and is designed to
take advantage of both global first-order information and local higher-order information. Indeed, when $\rho = 0$ and $\Lb_j$ is quadratic, \eqref{update} has the following closed form solution:
\begin{align*}
\xb^{t+1} = \xb^{t} - \nabla^2\Lb_1(\xb^t) ^{-1} \left(\textstyle\frac{1}{m}\sum_{j=1}^m \nabla \Lb_j(\xb^{t_j})\right),
\end{align*}
which is similar to a Newton updating step. The more general case has a \emph{proximal Newton} flavor; see, e.g., \cite{lee2014proximal} and the references therein. However, our method is different from their methods in the proximal term $\frac{\rho}{2}\|\xb - \xb^t\|^2$ as well as the first order term. Intuitively, if we have a first-order approximation
\begin{align} 
\Lb_1(\xb) \approx \Lb_1(\xb^t) + \left\langle\nabla\Lb_1(\xb^{t_1}), \xb - \xb^t\right\rangle, 
\end{align}
then $\eqref{update}$ reduces to
\begin{align}\label{updatered}
\nonumber \textstyle \xb^{t+1} &= \textstyle \argmin\limits_{\xb} \Big\langle\frac{1}{m}\sum_{j \in [m]} \nabla \Lb_j(\xb^{t_j}), \xb - \xb^t \Big\rangle \\
&+ h(\xb) + \frac{\rho}{2}\|\xb - \xb^t\|^2,  
\end{align}
which is essentially a first-order proximal gradient updating step.

We consider the convergence of the proposed approach under the asynchronous protocol where the master has the freedom to make updates with gradients from only a partial set of worker machines. 
We start with introducing important conditions that are used commonly in previous work \cite{hong2014distributed,chang2016asynchronous,hong2016convergence}. 

\begin{assumption}\label{Lipdiff}
	The function $\Lb_j(\xb)$ is differentiable and has Lipschitz continuous gradient for all $j\in [m]$, i.e.,
	\begin{align*}
	\|\nabla \Lb_j(\xb) - \nabla \Lb_j(\yb)\| \leq L\|\xb-\yb\|. 
	\end{align*}
\end{assumption}
The proof of the linear convergence relies on the following strong convexity assumption. 
\vspace{-0.1in}\begin{assumption}\label{strongconvex}
	For all $j\in [m]$, the function $\Lb_j$ is strongly convex with modulus $\sigma^2$, which means that
	\begin{align*}
	\Lb_j(x) > \Lb_j(y) + \langle \nabla \Lb_j(y),x - y\rangle + \frac{\sigma^2}{2} \|x-y\|^2,
	\end{align*}
	for all $\xb, \yb \in \RR^p$, $j \in[m]$. 
\end{assumption}

\begin{assumption}\label{convsubproblem}
	For all $t$, the parameter $\rho$ in \eqref{update} is chosen large enough such that:
	\begin{itemize}
		\item[I.] $\gamma(\rho) > 3L+ 2L\delta\tau$ and $\rho>\frac{2L\tau}{\delta}$, for some constant $\delta>0$,  where $\gamma(\rho)$ represents the convex modulus of the function $h(\xb) + \frac{\rho}{2}\|\xb-\xb^t\|^2$.  
		\item[II.] There exists a constant $\underline{L}$ such that
		\begin{align*}
		\Lb(\xb) > \underline{L}>-\infty  \hspace{1.1em} \forall \xb\in \RR^p. 
		\end{align*}
	\end{itemize}
\end{assumption}

Moreover the following concept is needed in the first part of Theorem~\ref{mainthm}. 

\begin{definition}
	We say a function $\cF(\xb)$ is coercive if 
	\begin{align*}
	\lim\limits_{\|\xb\| \rightarrow \infty} \cF(\xb) = + \infty. 
	\end{align*} 
\end{definition}

Define
\begin{align}
\textstyle \tilde{\nabla}_\xb\Lb\left(\xb^t\right) = \textstyle \xb^t - \textbf{Prox}_{h}\Big(\xb^t - \frac{1}{m} \sum\limits_{j\in [m]}\nabla \Lb_j(\xb^{t})\Big),
\end{align}
where $\textbf{Prox}_{h}$ is a proximal operator defined by $\textbf{Prox}_{h} [z]:= \argmin\limits_{x}h(x) + \frac{1}{2}\|x-z\|^2$. Usually $\tilde{\nabla}_\xb\Lb\left(\xb^t\right)$ is called the proximal gradient of $\Lb$; $\xb$ is a stationary point when $\tilde{\nabla}_\xb\Lb\left(\xb\right)=0$. 

Based on these assumptions, now we can present the main theorem. 
\begin{thmi}\label{mainthm}
	Suppose Assumption~\ref{boundeddelay}, \ref{Lipdiff}, and \ref{convsubproblem} are satisfied. Then we have the following claims for the sequence generated by Algorithm \ref{algor} (EDANNI). 
	\begin{itemize}
		\item[$\bullet$] {\bf (Boundedness of Sequence).}  The gap between $\xb^t$ and $\xb^{t+1}$ converges to $0$, i.e.,
		\begin{align*}
		\lim\limits_{t\rightarrow \infty}\xb^{t+1} - \xb^{t} = 0. 
		\end{align*}
		If $\Lb(\xb)$ is coercive, then the sequence $\{\xb^{t}\}$ generated by Algorithm \ref{algor} is bounded. 
		\item[$\bullet$] {\bf (Convergence to Stationary Points).} Every limit point of the iterates $\{\xb^t\}$ generated by Algorithm \ref{algor} is a stationary point of problem \eqref{origin}. Furthermore,  $\left\|\tilde{\nabla}_\xb\Lb(\xb^{t})\right\|\rightarrow 0$, as $t\rightarrow \infty$. 
		\item[$\bullet$] {\bf (Sublinear Convergence Rate).} Given $\epsilon>0$, let us define $T$ to be the first time for the optimality gap to reach below $\epsilon$, i.e.,
		\begin{align*}
		T:= \argmin_{t} \left\{ \left \|\tilde{\nabla}_\xb\Lb(\xb^{t}) \right\|< \epsilon \right \}.
		\end{align*}
		Then there exists a constant $\nu>0$ such that  
		\begin{align*}
		T \leq \frac{\nu}{\epsilon} + 1, 
		\end{align*}
		where $\nu$ equals to a positive constant times $\left(2(2+\rho)^2 +8L^2\tau \right)/\min\left\{\frac{\gamma(\rho)}{2}- \frac{3L}{2} - L\delta\tau, \frac{\rho}{2}- \frac{L\tau}{\delta} \right\}$ for some $\delta>0$. Therefore, the optimality gap $\left\|\tilde{\nabla}_\xb\Lb(\xb^{t})\right\|$ converges to $0$ in a sublinear manner. 
	\end{itemize}
	
\end{thmi}

\begin{remark}
	The theorem suggests that the iterates $\{\xb^t\}$ may or may not be bounded without the coerciveness property of $\Lb(\xb)$. However, it guarantees that the optimality measure $\left\|\tilde{\nabla}_x\Lb(\xb^{t})\right\|$ converges to $0$ sublinearly. We remark that \cite{li2014communication} also analyzed the convergence of a proximal gradient method based communication-efficient algorithm for nonconvex problems, but they did not give a specific convergence rate. Note that such sublinear complexity bound is tight when applying first-order methods for nonconvex unconstrained problems (see \cite{ghadimi2016mini,cartis2010complexity}). 
\end{remark}
Let us define
\begin{align*}
{\textstyle \text{F}(\xb,\xb^t) :=  \frac{1}{m}\sum_{j\in [m]}\Lb_j(\xb) + \frac{\rho}{2}\left\| \xb- \xb^t\right\|^2 + h(\xb). }
\end{align*}
The gap between $\xb^t$ and $\xb^{t+1}$ is denoted by $\Delta^{(t)} := \xb^{t+1} - \xb^{t}$, for $t \in \NN $. The proof of Theorem~\ref{mainthm} relies on Lemma~\ref{partdescent}, \ref{descentlemma}, and \ref{Flowbound} in the following.  

\begin{lemmai}\label{partdescent}
	Suppose Assumption \ref{Lipdiff} and Assumption~\ref{convsubproblem} (I) are satisfied. then the following is true for iterates $\{\xb^t\}$ generated by Algorithm \ref{algor} (EDANNI)
	\begin{align}\label{insertinq}
	\nonumber &\frac{\rho}{2}\|\xb^{t+1}-\xb^{t}\|^2 + h(\xb^{t+1}) - \frac{\rho}{2}\|\xb^{t}-\xb^{t-1}\|^2 - h(\xb^{t}) \\
	\nonumber&\leq -\Bigg \langle \nabla\Lb_1(\xb^{t+1}) + \frac{1}{m}\sum_{j=1}^m \nabla\Lb_j(\xb^{t_j})- \nabla\Lb_1(\xb^{t_1}), \\
	& \Delta^{(t)} \Bigg \rangle - \frac{\gamma(\rho)}{2}\big\|\Delta^{(t)}\big\|^2 -  \frac{\rho}{2}\big\|\Delta^{(t-1)}\big\|^2.
	\end{align}
\end{lemmai}

\begin{lemmai}	\label{descentlemma}
	Under the assumptions of Theorem \ref{mainthm} for any $\delta>0$ we have 
	
	$\text{F}(\xb^{t+1}, \xb^{t}) - \text{F}(\xb^{t}, \xb^{t-1}) $
	\begin{align} \label{adjacdescent1}
	\nonumber\leq \left(\frac{3L}{2} - \frac{\gamma(\rho)}{2}+L\delta\tau\right)\big\|\Delta^{(t)}\big\|^2 - \frac{\rho}{2}\big\|\Delta^{(t-1)}\big\|^2 \\
	+ \frac{L}{\delta}\sum_{k = 1}^{\tau} \big\|\Delta^{(t-k)}\big\|^2. 
	\end{align}
\end{lemmai}

\begin{lemmai}\label{Flowbound}
	Suppose Assumption \ref{convsubproblem} is satisfied. Then for $\xb^t$ generated by (EDANNI), there exists some constants $\underline{\Ft}$ and $\bar{\Ft}$ such that
	\begin{align*}
	\hspace{1.1em}  +\infty > \bar{\Ft}  > \Ft(\xb^{t+1}, \xb^{t}) > \underline{\Ft} > -\infty, \hspace{1.1em} \forall t \geq 0. 
	\end{align*}	
\end{lemmai}

The proofs of these lemmata are in the Appendix. Now in the following we prove Theorem \ref{mainthm}. 

\begin{proof}[\textbf{Proof of Theorem~\ref{mainthm}}]
We begin by establishing the first conclusion of the theorem. Summing inequality \eqref{adjacdescent1} in Lemma~\ref{descentlemma} over $t$ yields

$\text{F}(\xb^{T+1}, \xb^{T}) - \text{F}(\xb^{1}, \xb^{0})$
\begin{align*}
\leq \sum_{t=1}^T \left(\frac{3L}{2} - \frac{\gamma(\rho)}{2} + L\delta\tau\right)\big\| \Delta^{(t)}\big\|^2\\
+ \sum_{t=1}^T\left( \frac{L\tau}{\delta} -\frac{\rho}{2}\right)\big\| \Delta^{(t-1)}\big\|^2. 
\end{align*}
   
Now define 
\begin{align*}
c := \min\left\{\frac{\gamma(\rho)}{2}- \frac{3L}{2} - L\delta\tau, \frac{\rho}{2}- \frac{L\tau}{\delta} \right\}, 
\end{align*}
by Assumption~\ref{convsubproblem} we have $\gamma(\rho) > 3L+ 2L\delta\tau$ and $\rho>\frac{2L\tau}{\delta}$, therefore $c>0$. It holds that 
\begin{align}\label{Fdescenfinal}
\text{F}(\xb^{T+1}, \xb^{T}) - \text{F}(\xb^{1},\xb^{0}) \leq -c \sum_{t=0}^T \big\| \Delta^{(t)}\big\|^2. 
\end{align}

Note that by Lemma~\ref{Flowbound} the LHS of \eqref{Fdescenfinal} is bounded from below. By letting $T \rightarrow \infty$, it follows that
\begin{align*}
\big\|\Delta^{(t)}\big\| \rightarrow 0, \hspace{2em} t\rightarrow \infty. 
\end{align*}

Moreover, Lemma~\ref{Flowbound} shows that $\text{F}(\xb^{t+1},\xb^{t})$ is bounded, but due to the coerciveness assumption 
\begin{align}
\lim\limits_{\|\xb\|\rightarrow \infty} \frac{1}{m}\sum_{j \in [m]} \Lb_j(\xb)+h(\xb) + \frac{\rho}{2}\left\|\xb - \xb^t\right\|^2 = +\infty,
\end{align} 
so we know $\{\xb^{t+1}\}$ is bounded. Therefore the first conclusion is proved.

We now establish the second conclusion of the Theorem.  From \eqref{update}, we know that
\begin{align*}
\xb^{t+1} = \textbf{Prox}_{h} \Big[\xb^{t+1} - \Big( \nabla\Lb_1(\xb^{t+1}) + \frac{1}{m}\sum_{j \in [m]} \nabla \Lb_j(\xb^{t_j}) \\
- \nabla\Lb_1(\xb^{t_1}) +\rho \left( \xb^{t+1} - \xb^t\right) \Big)\Big],
\end{align*}
where $\textbf{Prox}_{h}$ is a proximal operator defined by $\textbf{Prox}_{h} [z]:= \argmin\limits_{x}h(x) + \frac{1}{2}\|x-z\|^2$.  This implies that
\begin{align}\label{proximalgradientbound} 
\nonumber&\Bigg\| \xb^t - \textbf{Prox}_{h}\Bigg(\xb^t - \frac{1}{m} \sum\limits_{j\in [m]}\nabla \Lb_j(\xb^{t})\Bigg)\Bigg\| \\
\nonumber& \leq \Bigg\|\xb^t - \xb^{t+1} + \xb^{t+1}\\
\nonumber& -  \textbf{Prox}_{h}\left(\xb^t - \frac{1}{m} \sum_{j\in [m]}\nabla \Lb_j(\xb^{t})\right) \Bigg\|\\ 
\nonumber&\leq  \left\|\xb^t - \xb^{t+1}\right\| \\
\nonumber& + \Bigg\| \textbf{Prox}_{h} \Bigg[\xb^{t+1} - \Bigg( \nabla\Lb_1(\xb^{t+1}) + \frac{1}{m}\sum_{j \in [m]} \nabla \Lb_j(\xb^{t_j}) \\
\nonumber& - \nabla\Lb_1(\xb^{t_1}) +\rho \big( \xb^{t+1} - \xb^t\big) \Bigg)\Bigg] \\
\nonumber& - \textbf{Prox}_{h}\left(\xb^t - \frac{1}{m} \sum_{j\in [m]}\nabla \Lb_j(\xb^{t})\right) \Bigg\|\\
\nonumber& \overset{(a)}{\leq}\Bigg\|(1+\rho)(\xb^{t+1} - \xb^t) + \frac{1}{m}\sum_{j \in [m]} \nabla \Lb_j(\xb^{t}) \\
\nonumber&   - \frac{1}{m}\sum_{j \in [m]} \nabla \Lb_j(\xb^{t_j}) - \left(\nabla \Lb_1(\xb^{t+1}) - \nabla \Lb_1(\xb^{t_1})\right) \Bigg\|\\
\nonumber&  +\big\|\Delta^{(t)}\big\| \\
\nonumber& \leq \left(2+\rho\right)\big\|\Delta^{(t)}\big\| + 2L\sum_{k=0}^\tau\big\|\Delta^{(t-k)}\big\| \\
& \longrightarrow 0, \hspace{2em} t\rightarrow \infty. 
\end{align} 
Note that here inequality $(a)$ holds because of the nonexpansiveness of the operator $\textbf{Prox}_{h}$. The last inequality follows from Assumption \ref{Lipdiff}.  

Let $\Xb^*$ be the set of stationary points of problem \eqref{origin}, and let 
\begin{align*}
\text{dist}\left(\xb^t, \Xb^*\right):=\min_{\hat{x}\in\Xb^*}{\|\xb^t-\hat{\xb}\|}
\end{align*}
denote the distance between $\xb^t$ and the set $\Xb^*$. Now we prove  
\begin{align*}
\lim\limits_{t\rightarrow\infty} \text{dist}\left(\xb^t, \Xb^*\right) = 0.
\end{align*}

Suppose there exists a subsequence $\{\xb^{t_k}\}$ of $\{\xb^{t}\}$ such that $\xb^{t_k}\rightarrow \hat{\xb}, ~ k\rightarrow\infty$ but 
\begin{align}\label{fanzheng}
\lim\limits_{k\rightarrow\infty}\text{dist}(\xb^{t_k},\Xb^*) \geq \gamma > 0.
\end{align}
Then it is obvious that $\lim\limits_{k\rightarrow\infty}\text{dist}(\xb^{t_k},\hat{\xb}) = 0$. Therefore there exists some $K(\gamma)>0$, such that
\begin{align}\label{distdef}
\|\xb^{t_k} - \hat{\xb}\| \leq \frac{\gamma}{2}, \hspace{1.1em} k > K(\gamma). 
\end{align}

On the other hand, from \eqref{proximalgradientbound} and the lower semi-continuity of $h(\xb)$ we have $\hat{\xb} \in \Xb^*$, so by the definition of the distance function we have
\begin{align} \label{distineq}
\text{dist}(\xb^{t_k},\Xb^*) \leq \text{dist}(\xb^{t_k},\hat{\xb}). 
\end{align}
Combining \eqref{distdef} and \eqref{distineq}, we must have 
\begin{align*}
\text{dist}(\xb^{t_k},\Xb^*) \leq \frac{\gamma}{2},\hspace{1.1em} k>K(\gamma).
\end{align*}
This contradicts to \eqref{fanzheng}, so the second result is proved.

We finally prove the third conclusion of the Theorem. Summing \eqref{proximalgradientbound} over $t$ yields
\begin{align}\label{gradbound}
\nonumber &\sum\limits_{t = 0}^T \Bigg\| \xb^t - \textbf{Prox}_{h}\left(\xb^t - \frac{1}{m} \sum\limits_{j\in [m]}\nabla \Lb_j(\xb^{t})\right)\Bigg\|^2 \\
\nonumber & \leq \sum_{t = 0}^T 2(2+\rho)^2\big\|\Delta^{(t)}\big\|^2 + 2(2L)^2 \sum_{t = 0}^T \sum_{k = 0}^\tau \big\|\Delta^{(t-k)}\big\|^2\\
& \leq \left(2(2+\rho)^2 +8L^2\tau \right) \sum_{t=0}^T \big\|\Delta^{(t)}\big\|^2. 
\end{align}
Combining \eqref{Fdescenfinal} and \eqref{gradbound} we have 
\begin{align*}
\sum\limits_{t = 0}^T\left\|\tilde{\nabla}\Lb\left(\xb^t\right)\right\|^2 \leq \frac{\mu}{c}\left(\text{F}(\xb^1,\xb^0) - \text{F}(\xb^{T+1},\xb^T)\right),
\end{align*}
where $\mu:=\left(2(2+\rho)^2 +8L^2\tau \right)$. 

Let $T(\epsilon) := \min\left\{t\mid\left\|\tilde{\nabla}\Lb(\xb^t)\right\|\leq \epsilon, t\geq 0\right\}$. Then the above inequality implies 
\begin{align*}
T(\epsilon) \epsilon \leq \frac{\mu}{c}\left(\text{F}(\xb^1,\xb^0) - \text{F}(\xb^{T+1},\xb^T)\right).
\end{align*}
Thus it follows that
\begin{align*}
\epsilon \leq \frac{C\cdot \left(\text{F}(\xb^1,\xb^0) - \underline{\text{F}}\right)}{T(\epsilon)},
\end{align*}
where $C:= \frac{\mu}{c}>0$, proving Theorem \ref{mainthm}. 
\end{proof}

Besides the convergence in the nonconvex setting, in the following theorem we show that the proposed algorithm converges linearly if $\Lb_j$ is strongly convex. Quite interestingly, comparing with the results of \cite{pmlr-v70-wang17f,jordan2016communication,ren2017communication}, here the linear convergence is established without any statistical assumption of the data stored in each local machine. 
\begin{thmi}\label{mainthm1}
	Suppose Assumption~\ref{boundeddelay}, \ref{Lipdiff}, and \ref{strongconvex} are satisfied. If $\rho$ is sufficiently large such that 
	\begin{align*}
	\frac{\delta_1 L + \frac{\rho}{2}(1+\delta_1)}{\frac{\rho}{2}(1+\delta_1)+\delta_1} +\frac{3L}{2}-\frac{\rho}{2}+L\delta\tau < 0
	\end{align*}
	and
	\begin{align*}
	& \frac{\delta_1 L + \frac{\rho}{2}(1+\delta_1)}{\frac{\rho}{2}(1+\delta_1)+\delta_1}
	+\frac{3L}{2}-\frac{\rho}{2}+L\delta\tau\\
	& -\frac{\rho \eta}{2}+\Big(\frac{L}{\delta}+ \frac{\frac{\delta_1}{2}L^2\tau}{\frac{\rho}{2}(1+\delta_1) + \delta_1 }\Big)\frac{\eta^\tau-1}{\eta-1} <0,
	\end{align*}	
	for some $\delta>0$ and $\delta_1>(2L+\rho+1)/\sigma^2$, then it holds for the sequence generated by (EDANNI) that
	\begin{align*}
	0 &\leq \Ft(\xb^{t+1},\xb^{t}) - \Ft(\xb^*,\xb^*) \\
	& \leq \frac{1}{\eta^{t}} (\Ft(\xb^{1},\xb^{0}) - \Ft(\xb^*,\xb^*) ),
	\end{align*} 
	where $\eta := 1+\frac{1}{\frac{\rho}{2}(1+\delta_1)+\delta_1}$. 	
\end{thmi}\vspace{-0.5em}

Note the above conditions can be satisfied when $\rho$ is sufficiently larger than the order of $L$ and the exponential of $\tau$ and $\delta_1$ is larger than the order of $L/\sigma^2$. Theorem~\ref{mainthm1} asserts that with the strongly convexity of $\Lb_j$'s, the augmented optimality gap decreases linearly to zero under these conditions. Moreover, Assumption \ref{strongconvex} can be replaced by only requiring each $\Lb_j$  is convex and $\textstyle\frac{1}{m}\sum_{j\in [m]} \Lb_j$ is strongly convex with modulus $\sigma^2$. To prove Theorem~\ref{mainthm1}, we need the following lemma to bound the optimality gap of function $\Ft$. \vspace{-0.5em}

\begin{lemmai}\label{upbound_optm}
	Suppose Assumption~\ref{boundeddelay}, \ref{Lipdiff}, and \ref{strongconvex} hold and $\delta_1>(2L+\rho+1)/\sigma^2$ for some $\delta_1>0$, then it follows that
	\begin{align}\label{eqn18}
	\nonumber& \frac{1}{\frac{\rho}{2}(1+\delta_1)+\delta_1}\left( \Ft(\xb^{t+1}, \xb^t) - \Ft(\xb^*, \xb^*) \right)\\
	\nonumber &\leq \frac{\delta_1 L+ \frac{\rho}{2}(1+\delta_1)}{\frac{\rho}{2}(1+\delta_1)+\delta_1}\|\xb^t - \xb^{t+1} \|^2 \\
	& + \frac{1}{\frac{\rho}{2}(1+\delta_1)+\delta_1}\frac{\delta_1}{2m}L^2 \sum_{j\in[m]}\|\xb^{t_j}- \xb^{t}\|^2. \vspace{-0.5em}
	\end{align}	
\end{lemmai}
\vspace{-0.7em}

The proof of Lemma~\ref{upbound_optm} is in the Appendix. Now we begin to prove Theorem \ref{mainthm1}.
\vspace{-0.101in}

\begin{proof}[\textbf{Proof of Theorem~\ref{mainthm1}}]	
	We begin by defining $\tilde{\Delta}^{(t+1)} = \Ft(\xb^{t+1},\xb^{t}) - \Ft(\xb^*,\xb^*) $. Then from the proof of Lemma~\ref{descentlemma} it holds that	
	\begin{align}\label{bylemma}
	\nonumber\tilde{\Delta}^{(t+1)} \leq \tilde{\Delta}^{(t)} + \left(\frac{3L}{2}- \frac{\rho}{2}+ L\delta \tau\right) \big\|\Delta^{(t)}\big\|^2 \\
    - \frac{\rho}{2} \big\|\Delta^{(t-1)}\big\|^2 + \left(\frac{L}{\delta}\sum_{k=1}^\tau \big\|\Delta^{(t-k)}\big\|^2\right)
	\end{align}
	Note that from \eqref{eqn18} of Lemma~\ref{upbound_optm} we have 
	\begin{align}\label{eqn19}
	\nonumber \frac{1}{\frac{\rho}{2}(1+\delta_1)+\delta_1}\tilde{\Delta}^{(t+1)}\leq \frac{\delta_1 L+\frac{\rho}{2}(1+\delta_1)}{\frac{\rho}{2}(1+\delta_1)+\delta_1}\big\|\Delta^{(t)} \big\|^2 \\
    + \frac{1}{\frac{\rho}{2}(1+\delta_1)+\delta_1}\frac{\delta_1}{2m}L^2 \sum_{j\in[m]}\|\xb^{t_j}- \xb^{t}\|^2
	\end{align}	
	
	By combining \eqref{eqn19} and \eqref{bylemma}, we have the following bound of the LHS:  
	\begin{align}\label{decay_of_deltat}
	\nonumber&\left(1+\frac{1}{\frac{\rho}{2}(1+\delta_1)+\delta_1}\right)\tilde{\Delta}^{(t+1)} \\
	\nonumber&\leq \tilde{\Delta}^{(t)} + \Big[ \frac{\delta_1 L + \frac{\rho}{2}(1+\delta_1)}{\frac{\rho}{2}(1+\delta_1)+\delta_1}+\frac{3L}{2}-\frac{\rho}{2}+L\delta\tau \Big] \big\|\Delta^{(t)}\big\|^2\\
	\nonumber& -\frac{\rho}{2}\big\|\Delta^{(t-1)}\big\|^2+ \frac{L}{\delta}\sum_{k=1}^\tau \big\|\Delta^{(t-k)}\big\|^2\\
	& + \frac{1}{\frac{\rho}{2}(1+\delta_1)+\delta_1}\frac{\delta_1}{2m}L^2 \sum_{j\in[m]}\|\xb^{t_j}- \xb^{t}\|^2.
	\end{align} 
	Inequality \eqref{decay_of_deltat} gives us an relation between $\tilde{\Delta}^{(t+1)}$ and $\tilde{\Delta}^{(t)}$. 
    Let us define $\eta := 1+\frac{1}{\frac{\rho}{2}(1+\delta_1)+\delta_1}$ and  $$(P3) := \frac{\delta_1 L + \frac{\rho}{2}(1+\delta_1)}{\frac{\rho}{2}(1+\delta_1)+\delta_1}+\frac{3L}{2}-\frac{\rho}{2}+L\delta\tau,$$ then by applying \eqref{decay_of_deltat} recursively we have 
	\begin{align*}
	&\tilde{\Delta}^{(t+1)} \\
	&\leq \frac{1}{\eta}\tilde{\Delta}^{(t)} + \frac{1}{\eta}(P3)\big\|\Delta^{(t)}\big\|^2 - \frac{\rho}{2\eta}\big\|\Delta^{(t-1)}\big\|^2 \\
	&+ \frac{L}{\delta\eta}\sum_{k=1}^\tau\big\|\Delta^{(t-k)}\big\|^2 +\frac{1}{\eta} \frac{1}{\frac{\rho}{2}(1+\delta_1)+\delta_1}\frac{\delta_1L^2}{2m} \sum_{j\in[m]}\|\xb^{t_j}- \xb^{t}\|^2\\
	&\leq \frac{1}{\eta^2}\tilde{\Delta}^{(t-1)} + \frac{1}{\eta} \left(\frac{1}{\eta}(P3)\big\|\Delta^{(t-1)}\big\|^2 - \frac{\rho}{2\eta}\big\|\Delta^{(t-2)}\big\|^2 \right)\\
	&  +\frac{1}{\eta}(P3)\big\|\Delta^{(t)}\big\|-\frac{\rho}{2\eta}\big\|\Delta^{(t-1)}\big\|^2 \\
	& + \left(\frac{L}{\delta\eta}\sum_{k=1}^\tau\big\|\Delta^{(t-k)}\big\|^2 + \frac{L}{\delta\eta^2}\sum_{k=1}^\tau\big\|\Delta^{(t-1-k)}\big\|^2 \right)\\
	&  + \frac{1}{\eta^2}\frac{1}{\frac{\rho}{2}(1+\delta_1)+\delta_1}\frac{\delta_1}{2m}L^2  \sum_{l=0}^1 \frac{1}{\eta^{l+1}} \sum_{j\in[m]}\sum_{k=1}^\tau\big\|\Delta^{(t-l-k)}\big\|^2 \\
	& \cdots
	\end{align*}
	\begin{align*}
	&\leq  \frac{1}{\eta^{t}}\tilde{\Delta}^{(1)} + \frac{1}{\eta}(P3)\big\|\Delta^{(t)}\big\|^2 + \Big( \frac{1}{\eta^2}(P3) -\frac{\rho}{2\eta} \Big)\big\|\Delta^{(t-1)}\big\|^2 \\
	& + \big( \frac{1}{\eta^3}(P3) -\frac{\rho}{2\eta^2} \big)\big\|\Delta^{(t-2)}\big\|^2+ \cdots  \\
	& + \big( \frac{1}{\eta^{t+1}}(P3) -\frac{\rho}{2\eta^t} \big)\big\|\Delta^{(0)}\big\|^2 + \Big(\frac{L}{\delta\eta}\sum_{l=0}^t \frac{1}{\eta^l}\sum_{k=1}^\tau\big\|\Delta^{(t-l-k)}\big\|^2 \Big)\\
	&  + \frac{1}{\frac{\rho}{2}(1+\delta_1)+\delta_1}\frac{\delta_1}{2m}L^2 \sum_{l=0}^t \frac{1}{\eta^{l+1}} \sum_{j\in[m]}\sum_{k=1}^\tau\big\|\Delta^{(t-l-k)}\big\|^2\\
	&\leq \frac{1}{\eta^{t}}\tilde{\Delta}^{(1)} + \frac{1}{\eta}(P3)\big\|\Delta^{(t)}\big\|^2 + \left( \frac{1}{\eta^2}(P3) -\frac{\rho}{2\eta} \right)\big\|\Delta^{(t-1)}\big\|^2 \\
	& + \left( \frac{1}{\eta^3}(P3) -\frac{\rho}{2\eta^2} \right)\big\|\Delta^{(t-2)}\big\|^2+ \cdots \\
	& + \left( \frac{1}{\eta^{t+1}}(P3) -\frac{\rho}{2\eta^t} \right)\big\|\Delta^{(0)}\big\|^2 \\
	&+ \Big(\frac{L}{\delta}+ \frac{\frac{\delta_1}{2}L^2\tau}{\frac{\rho}{2}(1+\delta_1) + \delta_1 }\Big)\frac{\eta^\tau-1}{\eta-1}\sum_{l=1}^t\frac{1}{\eta^{l+1}}\big\|\Delta^{(t-l)}\big\|^2,
	\end{align*}\vspace{-0.3em}
	where we use the fact that		
	\begin{align*}	
	&\sum_{l=0}^t \frac{1}{\eta^{l+1}}\sum_{k=1}^\tau\big\|\Delta^{(t-l-k)}\big\|^2 \\
	&\hspace{3em}=\eta^{-t-1}\sum_{l=0}^t \frac{\eta^t}{\eta^l}\sum_{k=1}^\tau\big\|\Delta^{(t-l-k)}\big\|^2\\
	&\hspace{3em}= \eta^{-t-1}\sum_{j=0}^t\eta^j\sum_{k=1}^\tau\big\|\Delta^{(j-k)}\big\|^2\\
	&\hspace{3em} \overset{(h)}{\leq} \eta^{-t-1}\sum_{j=0}^{t-1}\eta^{j+1}(1+\eta+\cdots+\eta^{\tau-1})\big\|\Delta^{(j)}\big\|^2\\
	&\hspace{3em} \leq  \frac{\eta^\tau-1}{\eta-1}\sum_{j=0}^{t-1}\frac{1}{\eta^{t-j}}\big\|\Delta^{(j)}\big\|^2\\
	&\hspace{3em} \leq \frac{\eta^\tau-1}{\eta-1}\sum_{l=1}^{t}\frac{1}{\eta^{l+1}}\big\|\Delta^{(t-l)}\big\|^2.
	\end{align*}
	The inequality $(h)$ holds because the coefficient of $\big\|\Delta^{(j)}\big\|^2$ in the summation is less than $\eta^{j+1}(1+\eta+\cdots+\eta^{\tau-1})$.
	
	Therefore if $\rho>0$ satisfies that
	\begin{align}
	(P3) := \frac{\delta_1 L + \frac{\rho}{2}(1+\delta_1)}{\frac{\rho}{2}(1+\delta_1)+\delta_1} +\frac{3L}{2}-\frac{\rho}{2}+L\delta\tau < 0,
	\end{align}
	and
	\begin{align} 
	\nonumber&(P3) -\frac{\rho \eta}{2} + \Big(\frac{L}{\delta}+ \frac{\frac{\delta_1}{2}L^2\tau}{\frac{\rho}{2}(1+\delta_1) + \delta_1 }\Big)\frac{\eta^\tau-1}{\eta-1} \\
	\nonumber &= \frac{\delta_1 L + \frac{\rho}{2}(1+\delta_1)}{\frac{\rho}{2}(1+\delta_1)+\delta_1}
	+\frac{3L}{2}-\frac{\rho}{2}+L\delta\tau\\
	& -\frac{\rho \eta}{2}+\Big(\frac{L}{\delta}+ \frac{\frac{\delta_1}{2}L^2\tau}{\frac{\rho}{2}(1+\delta_1) + \delta_1 }\Big)\frac{\eta^\tau-1}{\eta-1} <0, 
	\end{align}
	then we have 
	\begin{align*}
	0 \leq \tilde{\Delta}^{(t+1)}\leq \frac{1}{\eta^{t}}\tilde{\Delta}^{(1)}. 
	\end{align*}
	The conclusion is proved. 
\end{proof}

\vspace{-0.3in}
\subsection{Inexactly Solving the Subproblems}

In this section we discuss the case where subproblem~\eqref{update} is not solved exactly. The motivation is that in some practical applications, it may not be easy to exactly minimize the objective function. The following analysis shows that the convergence still holds true when there are small errors in solving the subproblems, thus implying the robustness of the proposed algorithm. Specifically, we assume subproblem~\eqref{update} is solved with some error at iteration $t$; that is, there is an error $\epsilon^t$ such that 
\begin{align}\label{inexact_update1}
\nonumber \textstyle \epsilon^t \in \nabla\Lb_1(\xb^{t+1}) + \frac{1}{m}\sum_{j \in [m]} \nabla \Lb_j(\xb^{t_j})- \nabla\Lb_1(\xb^{t_1}) \\
+\partial h(\xb^{t+1}) +\rho \left( \xb^{t+1} - \xb^t\right),
\end{align}
which is equivalent to 
\begin{align}\label{inexact_update}
\nonumber \textstyle \xb^{t+1} &= \textbf{Prox}_{h} \Big[\xb^{t+1} - \Big( \nabla\Lb_1(\xb^{t+1}) + \frac{1}{m}\sum_{j \in [m]} \nabla \Lb_j(\xb^{t_j}) \\
& - \nabla\Lb_1(\xb^{t_1}) +\rho \left( \xb^{t+1} - \xb^t\right) - \epsilon^t \Big)\Big].
\end{align}

First, we introduce the following assumption that gives the bound of the error term. 
\begin{assumption}\label{errorbound}
The error term in \eqref{inexact_update1} satisfies 
\begin{align*}
\|\epsilon^t\|^2 < c_1 \big\|\Delta^{(t-1)}\big\|^2, ~~\text{for}~~t>0.
\end{align*}
\end{assumption}
This assumption requires that the error in solving the subproblem is bounded by a constant times the progress of $\xb^t$ in the previous iteration. Note that when $\Delta^{(t-1)}:=\xb^t-\xb^{t-1} = 0$, it holds that $\xb^t$ is a stationary point in the nonconvex scenario and $\xb^t = \xb^*$ in the strongly convex scenario. 
Following the proof steps in Lemma~\ref{partdescent}, we have
\begin{align*}
\nonumber &\frac{\rho}{2}\|\xb^{t+1}-\xb^{t}\|^2 + h(\xb^{t+1}) - \frac{\rho}{2}\|\xb^{t}-\xb^{t-1}\|^2 - h(\xb^{t}) \\
\nonumber&\leq -\Big \langle \nabla\Lb_1(\xb^{t+1}) + \frac{1}{m}\sum_{j=1}^m \nabla\Lb_j(\xb^{t_j})- \nabla\Lb_1(\xb^{t_1}) - \epsilon^t, \\
& \Delta^{(t)} \Big \rangle - \frac{\gamma(\rho)}{2}\big\|\Delta^{(t)}\big\|^2 -  \frac{\rho}{2}\big\|\Delta^{(t-1)}\big\|^2.
\end{align*}

In the second step, for the descent of function $\Ft$ similar to Lemma~\ref{descentlemma} it holds that for any $\delta>0$ 

$\text{F}(\xb^{t+1}, \xb^{t}) - \text{F}(\xb^{t}, \xb^{t-1}) $
\begin{align}\label{inexact_descentF} 
\textstyle\nonumber\leq \left(\frac{3L}{2} - \frac{\gamma(\rho)}{2}+L\delta\tau\right)\big\|\Delta^{(t)}\big\|^2 - \frac{\rho}{2}\big\|\Delta^{(t-1)}\big\|^2 \\
\textstyle+ \frac{L}{\delta}\sum_{k = 1}^{\tau} \big\|\Delta^{(t-k)}\big\|^2 + \langle\epsilon^t,\Delta^{(t)}\rangle.
\end{align}
Therefore Lemma~\ref{Flowbound} still holds true by Assumption~\ref{errorbound} and \eqref{inexact_descentF}. Now the first conclusion of Theorem~\ref{mainthm} can be proved. Summing up inequality \eqref{inexact_descentF} over $t$ yields
\vspace{0.3em}

$\text{F}(\xb^{T+1}, \xb^{T}) - \text{F}(\xb^{1}, \xb^{0})$ \vspace{-0.2em}
\begin{align*}
& \leq \sum_{t=1}^T \left(\frac{3L}{2} - \frac{\gamma(\rho)}{2} + L\delta\tau\right)\big\| \Delta^{(t)}\big\|^2\\
&+ \sum_{t=1}^T\left( \frac{L\tau}{\delta} -\frac{\rho}{2}\right)\big\| \Delta^{(t-1)}\big\|^2 + \sum_{t=1}^T \langle\epsilon^t, \Delta^{(t)}\rangle\\
& \leq \sum_{t=1}^T \left(\frac{3L}{2} - \frac{\gamma(\rho)}{2} + L\delta\tau\right)\big\| \Delta^{(t)}\big\|^2\\
&+ \sum_{t=1}^T\left( \frac{L\tau}{\delta} -\frac{\rho}{2}\right)\big\| \Delta^{(t-1)}\big\|^2 + \frac{1}{2}\sum_{t=1}^T (\|\epsilon^t\|^2 + \big\|\Delta^{(t)} \big\|^2)\\
&\leq \sum_{t=1}^T \left(\frac{3L}{2} - \frac{\gamma(\rho)}{2} + L\delta\tau +\frac{1}{2}\right)\big\| \Delta^{(t)}\big\|^2\\
&+ \sum_{t=1}^T\left( \frac{L\tau}{\delta} -\frac{\rho}{2} + \frac{1}{2}c_1\right)\big\| \Delta^{(t-1)}\big\|^2,
\end{align*} 
where in the last inequality we use Assumption~\ref{errorbound}. 

Now define $\tilde{c} := \min \big\{\frac{\gamma(\rho)}{2}- \frac{3L}{2} - L\delta\tau -\frac{1}{2}, \frac{\rho}{2}- \frac{L\tau}{\delta}-\frac{1}{2}c_1 \big\}$.
Assume that \vspace{-0.5em}
\begin{align}\label{newrho_bd}
\gamma(\rho) > 3L+ 2L\delta\tau+1 \hspace{1em} \text{and} \hspace{1em} \rho>\frac{2L\tau}{\delta}+c_1,
\end{align} 
then we have $\tilde{c}>0$. Therefore
\begin{align}\label{Fdescenfinal_inexact}
\text{F}(\xb^{T+1}, \xb^{T}) - \text{F}(\xb^{1},\xb^{0}) \leq -\tilde{c} \sum_{t=0}^T \big\| \Delta^{(t)}\big\|^2. 
\end{align}
Note that by Lemma~\ref{Flowbound} the LHS of \eqref{Fdescenfinal_inexact} is bounded from below. It follows that
\begin{align*}
\big\|\Delta^{(t)}\big\| \rightarrow 0, \hspace{2em} t\rightarrow \infty. 
\end{align*}

We now establish the second conclusion of Theorem~\ref{mainthm}. From \eqref{inexact_update} we know that 
\begin{align*}
\xb^{t+1} = \textbf{Prox}_{h} \Big[\xb^{t+1} - \Big( \nabla\Lb_1(\xb^{t+1}) + \frac{1}{m}\sum_{j \in [m]} \nabla \Lb_j(\xb^{t_j}) \\
- \nabla\Lb_1(\xb^{t_1}) +\rho \left( \xb^{t+1} - \xb^t\right) +\epsilon^t \Big)\Big],
\end{align*}
where $\textbf{Prox}_{h}$ is a proximal operator defined by $\textbf{Prox}_{h} [z]:= \argmin\limits_{x}h(x) + \frac{1}{2}\|x-z\|^2$.  This implies that
\begin{align}\label{proximalgradientbound1}
\nonumber& \Bigg\| \xb^t - \textbf{Prox}_{h}\left(\xb^t - \frac{1}{m} \sum\limits_{j\in [m]}\nabla \Lb_j(\xb^{t})\right)\Bigg\|\\
\nonumber& \leq \Bigg\|\xb^t - \xb^{t+1} + \xb^{t+1}\\
\nonumber& -  \textbf{Prox}_{h}\left(\xb^t - \frac{1}{m} \sum_{j\in [m]}\nabla \Lb_j(\xb^{t})\right) \Bigg\|\\ 
\nonumber&\leq  \left\|\xb^t - \xb^{t+1}\right\| \\
\nonumber& + \Bigg\| \textbf{Prox}_{h} \Bigg[\xb^{t+1} - \Bigg( \nabla\Lb_1(\xb^{t+1}) + \frac{1}{m}\sum_{j \in [m]} \nabla \Lb_j(\xb^{t_j}) \\
\nonumber& - \nabla\Lb_1(\xb^{t_1}) +\rho \big( \xb^{t+1} - \xb^t\big) -\epsilon^t\Bigg)\Bigg] \\
\nonumber& - \textbf{Prox}_{h}\left(\xb^t - \frac{1}{m} \sum_{j\in [m]}\nabla \Lb_j(\xb^{t})\right) \Bigg\|\\
\nonumber&\overset{(\tilde{a})}{\leq}\Bigg\|(1+\rho)(\xb^{t+1} - \xb^t) -\epsilon^t + \frac{1}{m}\sum_{j \in [m]} \nabla \Lb_j(\xb^{t}) \\
\nonumber&   - \frac{1}{m}\sum_{j \in [m]} \nabla \Lb_j(\xb^{t_j}) - \left(\nabla \Lb_1(\xb^{t+1}) - \nabla \Lb_1(\xb^{t_1})\right) \Bigg\|\\
\nonumber&  +\big\|\Delta^{(t)}\big\| \\
\nonumber& \leq \left(2+\rho\right)\big\|\Delta^{(t)}\big\| + 2L\sum_{k=0}^\tau\big\|\Delta^{(t-k)}\big\| + c_1^{\frac{1}{2}}\big\|\Delta^{(t-1)}\big\| \\
& \longrightarrow 0, \hspace{2em} t\rightarrow \infty. 
\end{align} 
Note that here inequality $(\tilde{a})$ holds because of the nonexpansiveness of the operator $\textbf{Prox}_{h}$. The last inequality follows from Assumption \ref{Lipdiff}. Therefore as in the proof of Theorem~\ref{mainthm}, the second result holds. 

The rest analysis is the same as that of Theorem~\ref{mainthm}. Specifically it holds that 
\begin{align*}
\sum\limits_{t = 0}^T\left\|\tilde{\nabla}\Lb\left(\xb^t\right)\right\|^2 \leq \frac{\tilde{\mu}}{\tilde{c}}\left(\text{F}(\xb^1,\xb^0) - \text{F}(\xb^{T+1},\xb^T)\right),
\end{align*}
where $\tilde{\mu}:=3\left(2+\rho +2L\tau + c_1^{\frac{1}{2}} \right)$. 

Recall that $T(\epsilon) := \min\left\{t\mid\left\|\tilde{\nabla}\Lb(\xb^t)\right\|\leq \epsilon, t\geq 0\right\}$, thus it follows that
\begin{align*}
\epsilon \leq \frac{C\cdot \left(\text{F}(\xb^1,\xb^0) - \underline{\text{F}}\right)}{T(\epsilon)},
\end{align*}
where $C:= \frac{\tilde{\mu}}{\tilde{c}}>0$. Therefore we have the following corollary. 
\begin{cori}\label{corollary}
	Let Assumptions~\ref{boundeddelay} \ref{Lipdiff}, \ref{convsubproblem}(II), and \ref{errorbound} hold, and suppose $\rho$ satisfies \eqref{newrho_bd}. Then all conclusions in Theorem~\ref{mainthm} hold true for the sequence generated by (EDANNI) with subproblems being solved inexactly (as quantified above). 
\end{cori}

Results corresponding to Theorem~\ref{mainthm1} also holds in the scenario of solving the subproblems inexactly. Similar to Lemma~\ref{upbound_optm}, the optimality gap of function $F$ can be bounded by
\begin{align}\label{bd_inexact_strong_cvx}
	\nonumber& \frac{1}{\frac{\rho}{2}(1+\delta_1)+\delta_1}\left( \Ft(\xb^{t+1}, \xb^t) - \Ft(\xb^*, \xb^*) \right)\\
	\nonumber &\leq \frac{\delta_1 L}{\frac{\rho}{2}(1+\delta_1)+\delta_1}\|\xb^t-\xb^{t+1}\|^2 \\
	\nonumber&+ \frac{\frac{\rho}{2}(1+\delta_1)}{\frac{\rho}{2}(1+\delta_1)+\delta_1}\|\xb^t - \xb^{t+1} \|^2 \\
	\nonumber & + \frac{1}{\frac{\rho}{2}(1+\delta_1)+\delta_1}\frac{\delta_1}{2m}L^2 \sum_{j\in[m]}\|\xb^{t_j}- \xb^{t}\|^2\\
	& +\big(\frac{1}{2}\|\epsilon^t\|^2+ \frac{1}{2}\big\| \Delta^{(t)}\big\|^2\big)\frac{1}{\frac{\rho}{2}(1+\delta_1)+\delta_1}. 
\end{align}	
Following the steps of Lemma~\ref{descentlemma}, we have

$\text{F}(\xb^{t+1}, \xb^{t}) - \text{F}(\xb^{t}, \xb^{t-1}) $
\begin{align}\label{inexact_lemma2}
\nonumber\leq \left(\frac{3L}{2} - \frac{\gamma(\rho)}{2}+L\delta\tau\right)\big\|\Delta^{(t)}\big\|^2 - \frac{\rho}{2}\big\|\Delta^{(t-1)}\big\|^2 \\
+ \frac{L}{\delta}\sum_{k = 1}^{\tau} \big\|\Delta^{(t-k)}\big\|^2 + \langle\epsilon^t,\Delta^{(t)}\rangle.
\end{align}
Combining \eqref{bd_inexact_strong_cvx} and \eqref{inexact_lemma2} and then applying the Assumption~\ref{errorbound} leads to 
\begin{align*}
&\left(1+\frac{1}{\frac{\rho}{2}(1+\delta_1)+\delta_1}\right)\tilde{\Delta}^{(t+1)} \\
&\leq \tilde{\Delta}^{(t)} -\frac{\rho-c_1}{2}\big\|\Delta^{(t-1)}\big\|^2+ \frac{L}{\delta}\sum_{k=1}^\tau \big\|\Delta^{(t-k)}\big\|^2\\
& + \Big[ \frac{\delta_1 L + \frac{\rho}{2}(1+\delta_1)+\frac{1}{2}}{\frac{\rho}{2}(1+\delta_1)+\delta_1}+\frac{3L}{2}-\frac{\rho-1}{2} +L\delta\tau \Big] \big\|\Delta^{(t)}\big\|^2\\
& + \frac{1}{\frac{\rho}{2}(1+\delta_1)+\delta_1}\frac{\delta_1}{2m}L^2 \sum_{j\in[m]}\|\xb^{t_j}- \xb^{t}\|^2\\
&+\frac{c_1}{\frac{\rho}{2}(1+\delta_1)+\delta_1}\big\|\Delta^{(t-1)}\big\|^2.
\end{align*}

By a recursive argument similar to the proof of Theorem~\ref{mainthm1}, one can prove that if $\rho > 0$ satisfies that 
\begin{align}\label{para_bd_inexact1}
\frac{\delta_1 L + \frac{\rho}{2}(1+\delta_1)+\frac{1}{2}}{\frac{\rho}{2}(1+\delta_1)+\delta_1} +\frac{3L}{2}-\frac{\rho-1}{2}+L\delta\tau < 0
\end{align}
and
\begin{align}\label{para_bd_inexact2} 
\nonumber& \frac{\delta_1 L + \frac{\rho}{2}(1+\delta_1)+\frac{1}{2}}{\frac{\rho}{2}(1+\delta_1)+\delta_1}
+\frac{3L}{2}-\frac{\rho-1}{2}+L\delta\tau -\frac{(\rho-c_1) \eta}{2}\\
&  +\Bigg(\frac{L}{\delta}+ \frac{\frac{\delta_1}{2}L^2\tau}{\frac{\rho}{2}(1+\delta_1) + \delta_1 }\Bigg)\frac{\eta^\tau-1}{\eta-1} + \frac{c_1\eta}{\frac{\rho}{2}(1+\delta_1)+\delta_1}<0,
&\textstyle  
\end{align}
then we have
\begin{align*}
0 \leq \tilde{\Delta}^{(t+1)}\leq \frac{1}{\eta^{t}}\tilde{\Delta}^{(1)}. 
\end{align*}
In summary we have the following corollary. 
\begin{cori}\label{corollary1}
	Suppose Assumption~\ref{boundeddelay}, \ref{Lipdiff}, \ref{strongconvex}, and \ref{errorbound} are satisfied. If $\rho$ satisfies \eqref{para_bd_inexact1} and \eqref{para_bd_inexact2} for some $\delta>0$ and $\delta_1>(2L+\rho+1)/\sigma^2$, then for the sequence generated by (EDANNI) with subproblems being solved inexactly we have 
	\begin{align*}
	0 &\leq \Ft(\xb^{t+1},\xb^{t}) - \Ft(\xb^*,\xb^*) \\
	& \leq \frac{1}{\eta^{t}} \Big(\Ft(\xb^{1},\xb^{0}) - \Ft(\xb^*,\xb^*) \Big),
	\end{align*} 
	where $\eta := 1+\frac{1}{\frac{\rho}{2}(1+\delta_1)+\delta_1}$. 	
\end{cori}

\vspace{0.3em}
\section{Experiments}\label{experimentsection}

Now we test our algorithm on both a convex application (LASSO) and a nonconvex application (Sparse PCA). In both settings, we compare with various advanced algorithms: 
\begin{itemize}
	\item[(1)] Efficient Distributed Learning with the Parameter Server (Parameter Server): the state-of-the-art proximal gradient descent based framework with the parameter server proposed in \cite{li2014communication}.
	\item[(2)] Asynchronous Distributed ADMM (AD-ADMM): the ADMM based asynchronous algorithm proposed in \cite{chang2016asynchronous}. 
	\item[(3)] Efficient Distributed Algorithm for Nonconvex-Nonsmooth Inference (EDANNI): the proposed approach in this paper. 	
\end{itemize}
We first compare their communication cost, that is, the total number of transmissions between the master
and the workers. Then the effects of the asynchrony on the working time and the idle time are examined. 

\subsection{LASSO}

In this example, to demonstrate the convergence performance of the above algorithms in terms of communication rounds, we consider the following LASSO problem
\begin{align}
\argmin_{\wb \in \RR^p} \frac{1}{2mn}\sum_{j\in [m]} \sum_{i\in [n]}\left\| \xb_{ji}^T\wb - y_{ji}\right\|^2 + \theta \|\wb\|_1,
\end{align}
where $\theta>0$ is the coefficient of the regularizer. Note that from now on we switch the notation to let the unknown quantity be denoted by $\wb$ instead of $\xb$.   

The data $\{\xb_{ji}\}_{i\in [n], j\in[m]}$ is independently sampled from a multivariate Gaussian distribution with zero mean and covariance matrix $\Sigma$. For $r\in [p], t \in [p]$, the $rt$-th entry of covariance matrix is set to be: $|\Sigma_{rt}| = 0.5^{|r-t|}$. The corresponding  $y_{ji}$ is constructed by 
\begin{align*}
y_{ji} = \xb_{ji}^T\wb^* + \epsilon_{ji}, \hspace{1.1em} \forall j \in [m], i\in [n],
\end{align*}
where noise $\epsilon_{ji}$ is a zero mean Gaussian random variable with variance $0.01$. The true parameter $\wb^*$ is $s$-sparse where all the entries are zero except that the first
$s$ entries are i.i.d random variables from a uniform distribution in [0,1]. The sparsity $s$ is set to be $0.01\times p$, where $p$ is the dimension of data $\xb_{ji}$. 

Those algorithms are compared in the setting where $m= 20$, $n =500$, $p = 1000$, $s=10$, and $\theta = 0.01$. Even though Theorem~\ref{mainthm} suggests that $\rho$ should be a larger value, we find that $\rho=0$ works well for this case. Both the synchronous ~scenario ~and ~the asynchronous scenario are considered. To simulate the asynchronous case, in each iteration half of the workers are assumed to be arrived with probability $0.2$ and the other workers are assumed be arrived with probability $0.5$. Moreover, the maximum tolerable delay $\tau$ is set be $3$ and the master will update the variables once the workers $j$ with $d_j>\tau-1$ have arrived.

\begin{table} \renewcommand*{\arraystretch}{1}
	\centering
	\caption{\small Comparison of the communication cost for LASSO}
	{\scriptsize \begin{tabular}{|p{1cm}|c|c|c|c|}
			\hline
			\multirow{2}{*}{\textbf{Type}} & \multirow{2}{*}{{\vspace{0.5em}\textbf{Parameters}}} & \multirow{2}{*}{\textbf{Method}}& \multirow{2}{*}{\textbf{Communication} }\\
			&{(m,n,p,s,$\theta$)}&& \\\hline
			\multirow{9}{*}{\hspace{-0.7em} \textbf{{\tiny Synchronous}}} &	\multirow{3}{*}{(10,1000,500,5,0.01)}
			&\textbf{AD-ADMM}					&21.9\\\cline{3-4}
			&&\textbf{PS} &1.7\\\cline{3-4}
			&&\textbf{EDANNI}	  &\textbf{1}\\\cline{2-4}
			&	\multirow{3}{*}{(20,500,500,5,0.01)}
			&\textbf{AD-ADMM}				    	&34.1\\\cline{3-4}
			&&\textbf{PS} &1.3\\\cline{3-4}
			&&\textbf{EDANNI}	   &\textbf{1}\\\cline{2-4}
			&	\multirow{3}{*}{(20,500,1000,10,0.01)}
			&\textbf{AD-ADMM}					    &8.3\\\cline{3-4}
			&&\textbf{PS} &1.4\\\cline{3-4}
			&&\textbf{EDANNI}	   &\textbf{1}\\\hline
			\multirow{7}{*}{\hspace{-0.7em} \textbf{{\tiny Asynchronous}}} &	\multirow{3}{*}{(10,1000,500,5,0.01)}
			&\textbf{AD-ADMM}					&20.6\\\cline{3-4}
			&&\textbf{PS} &1.2\\\cline{3-4}
			&&\textbf{EDANNI}	   &\textbf{1}\\\cline{2-4}
			&	\multirow{3}{*}{(20,500,500,5,0.01)}
			&\textbf{AD-ADMM}				    	&35.1\\\cline{3-4}
			&&\textbf{PS} &1.3\\\cline{3-4}
			\multirow{1}{*}{\hspace{0.1em} \textbf{{\tiny($\tau = 10$)}}}&&\textbf{EDANNI}	   &\textbf{1}\\\cline{2-4}
			&	\multirow{3}{*}{(20,500,1000,10,0.01)}
			&\textbf{AD-ADMM}					    &6.7\\\cline{3-4}
			&&\textbf{PS} &1.5\\\cline{3-4}
			&&\textbf{EDANNI}	  &\textbf{1}\\\hline
	\end{tabular}} \label{table_lasso}
\end{table}

\vspace{-0.22in}
\begin{figure}[!tbh] 
	\centering 
	\vspace{0.01in}
	\begin{center}
		\hspace{-0.15in}\subfloat[{\scriptsize Synchronous, \textbf{(7.7, 1.3, 1)}}]{
			\includegraphics[width=0.46\linewidth]{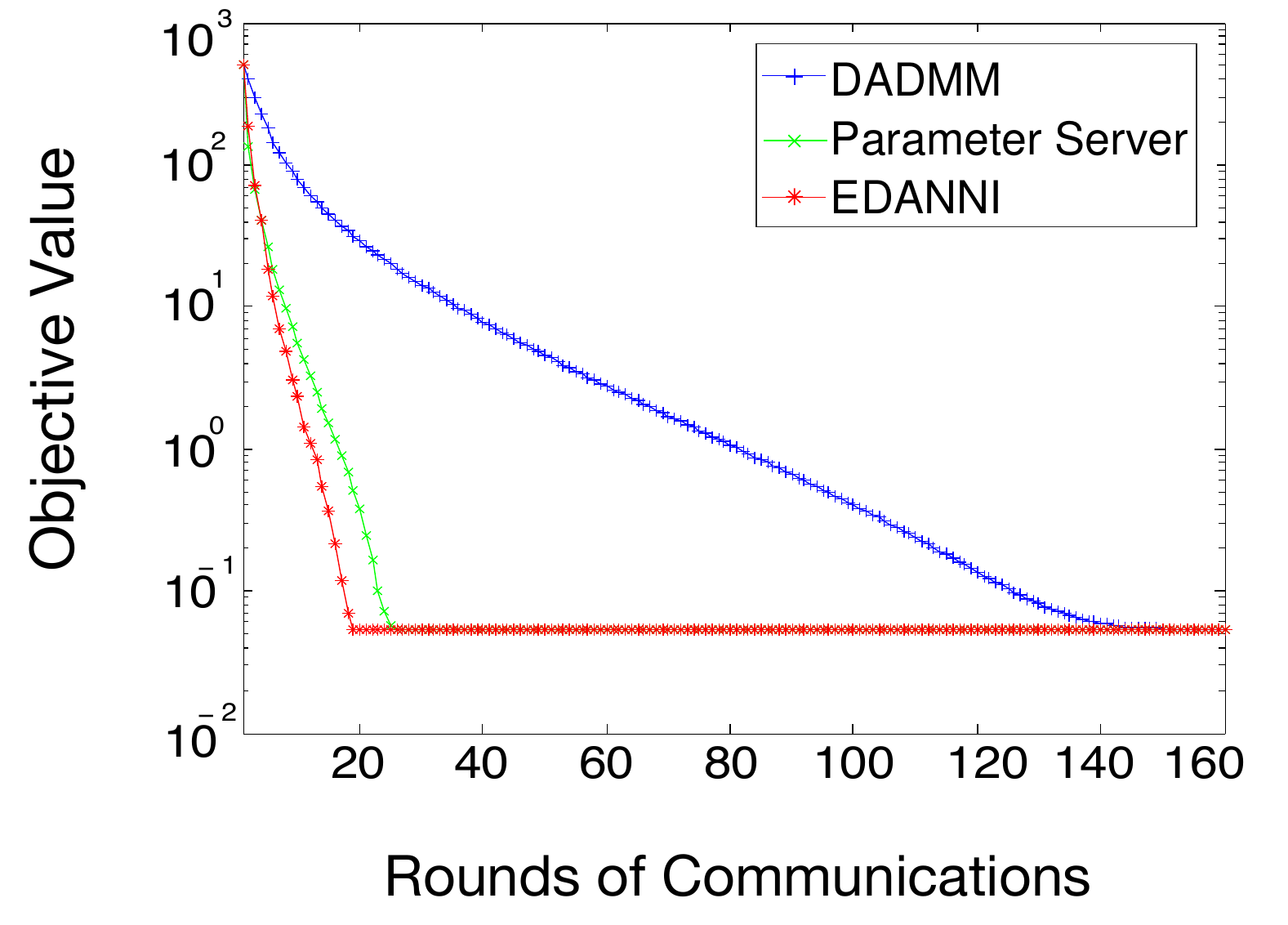}
		}\hspace{0.5em}
		\hspace{-0.11in}\subfloat[{\scriptsize Asynchronous, \textbf{(6.7, 1.5, 1)} }]{
			\includegraphics[width=0.46\linewidth]{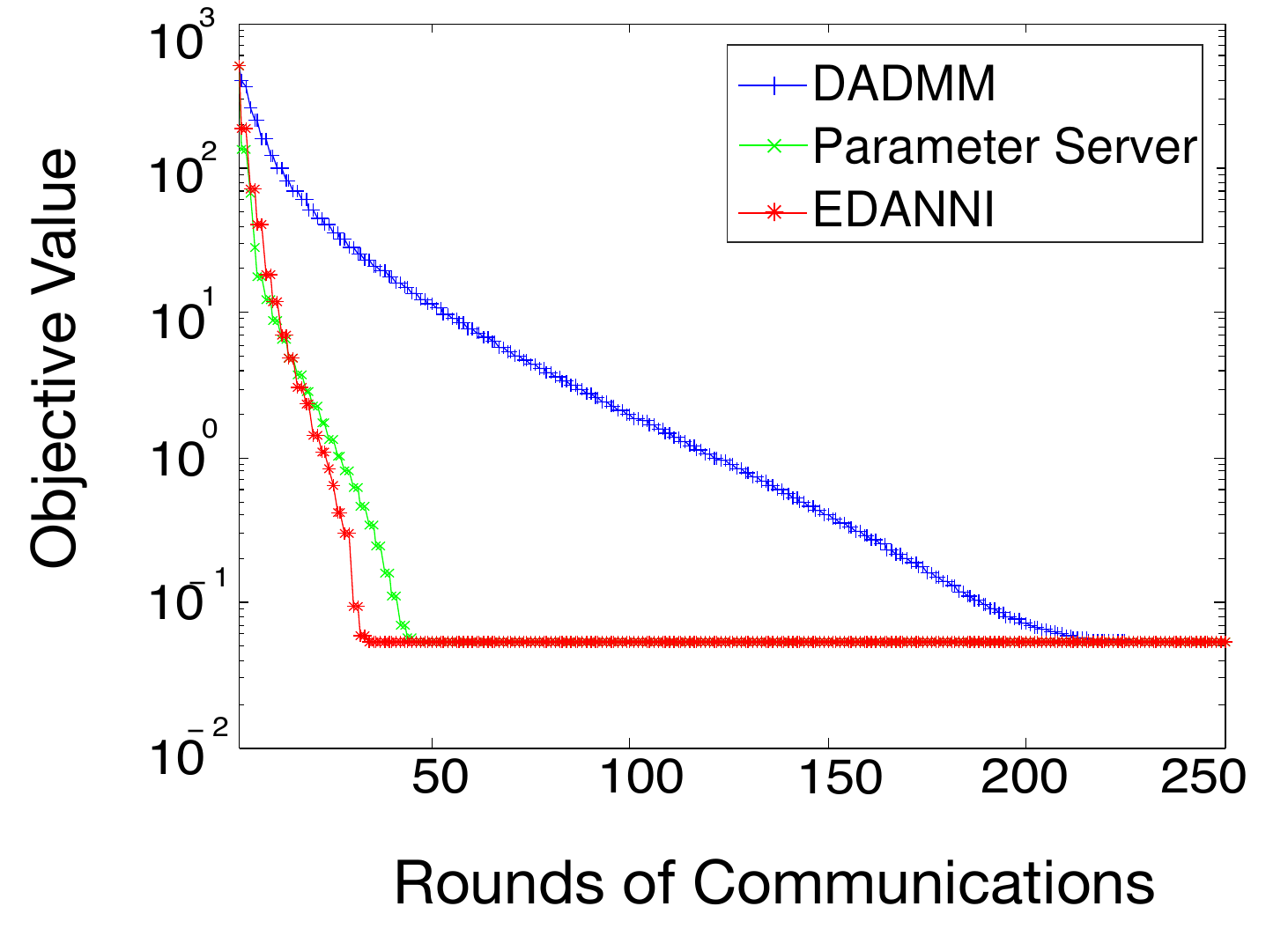}
		}
	\end{center}
	\vspace{-0.1in}	
	\caption{Comparison of candidate algorithms in LASSO when $m= 20$, $n =500$, $p = 1000$, $s=10$, and $\Xb \sim \cN(0,\Sigma)$. }\label{lassodiffalgrsyn}
	\vspace{-0.18in}
\end{figure} 
\begin{figure}[!tbh]
	\centering 
	\vspace{0.01in}
	\begin{center}
		\hspace{-0.07in}\subfloat[{\scriptsize Convergence for different $\tau$'s}]{
			\includegraphics[width=0.48\linewidth]{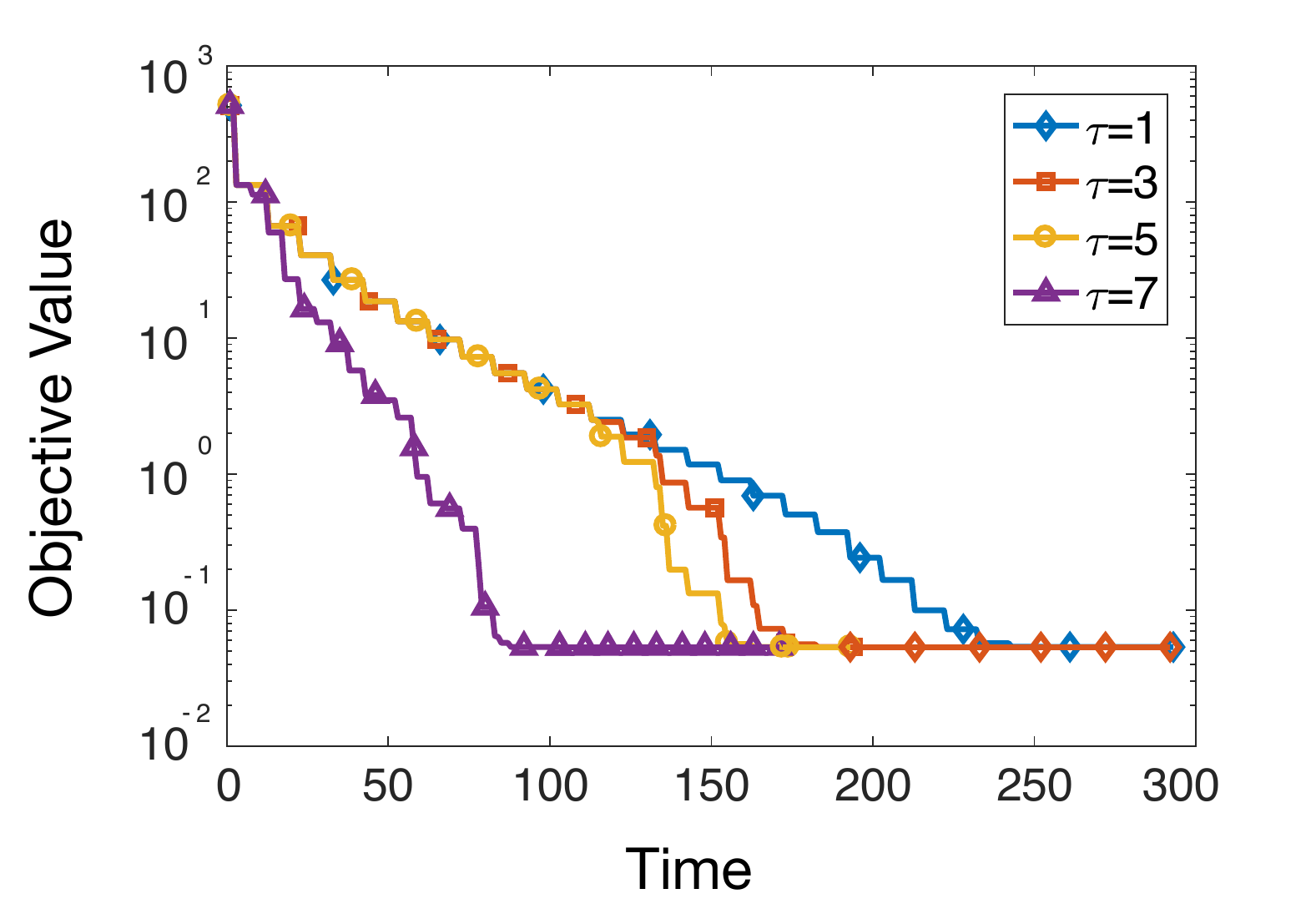}
		}\hspace{0.0em}
		\hspace{-0.08in}\subfloat[{\scriptsize Time table for different  $\tau$'s}]{
			\includegraphics[width=0.48\linewidth]{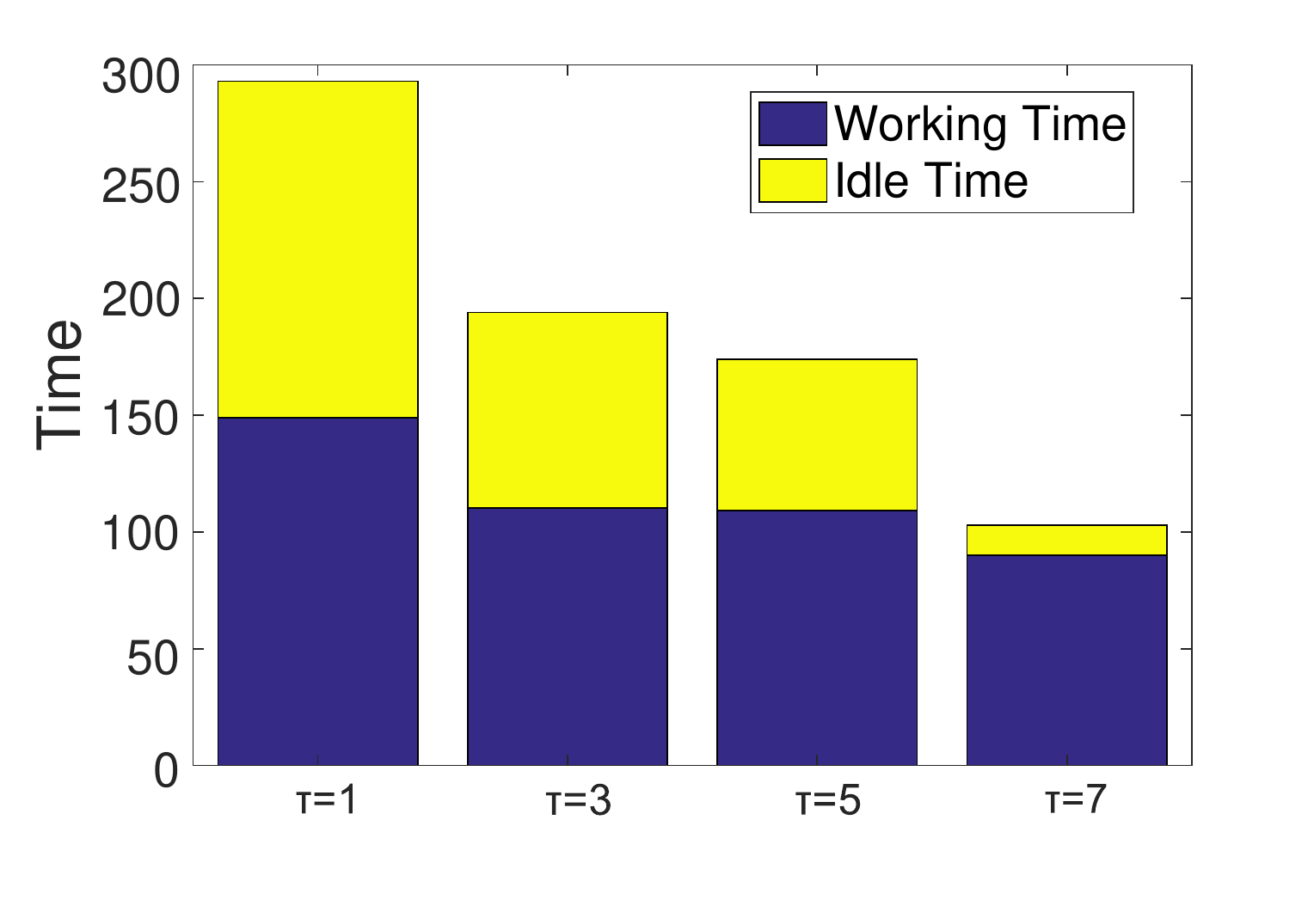}
		} 
	\end{center} 
	\vspace{-0.1in}	
	\caption{Comparison of EDANNI in LASSO with different settings of $\tau$ when $\rho=0$. }\label{difftau}
	\vspace{-0.1in}
\end{figure}
 

One can observe from Figure~\ref{lassodiffalgrsyn} that the proposed algorithm indeed converges faster than AD-ADMM and Parameter Server in terms of communication rounds,  in both the synchronous scenario and the asynchronous scenario. The triple of numbers in each figure's caption indicates the communication cost needed for AD-ADMM, Parameter Server, and EDANNI to attain the minimum objective value with error less than $10^{-6}$. For simplicity, we scale the communication complexity of EDANNI to 1. The results for other settings of $m, n, p, s$ are summarized in Table~\ref{table_lasso}. It is shown in these results that EDANNI is the most communication-efficient among the three algorithms.

Figure~\ref{difftau} shows the performance of the proposed approach when we choose different maximum tolerable delay $\tau$. It can be observed from Figure~\ref{difftau} (a) that the convergence rate varies much with different values of $\tau$. Basically, EDANNI converges faster when $\tau$ is larger, and converges relatively slower when $\tau = 0$ (i.e., the synchronous case). The results in Figure~\ref{difftau} (b) shows that the ratio of the computing time over the idle time increases when the delay bound $\tau$ becomes larger, therefore speeding up the convergence. Here the distributed implementation is simulated on a single machine by randomly setting the computation speed for each node from a uniform$[1,10]$ distribution.

\subsection{Sparse PCA}
\begin{figure}[!tbh]
	\centering 
	\vspace{0.01in}
	\begin{center}
		\hspace{-0.125in}\subfloat[{\scriptsize Synchronous, \textbf{(2.5, 2.0, 1)}}]{
			\includegraphics[width=0.46\linewidth]{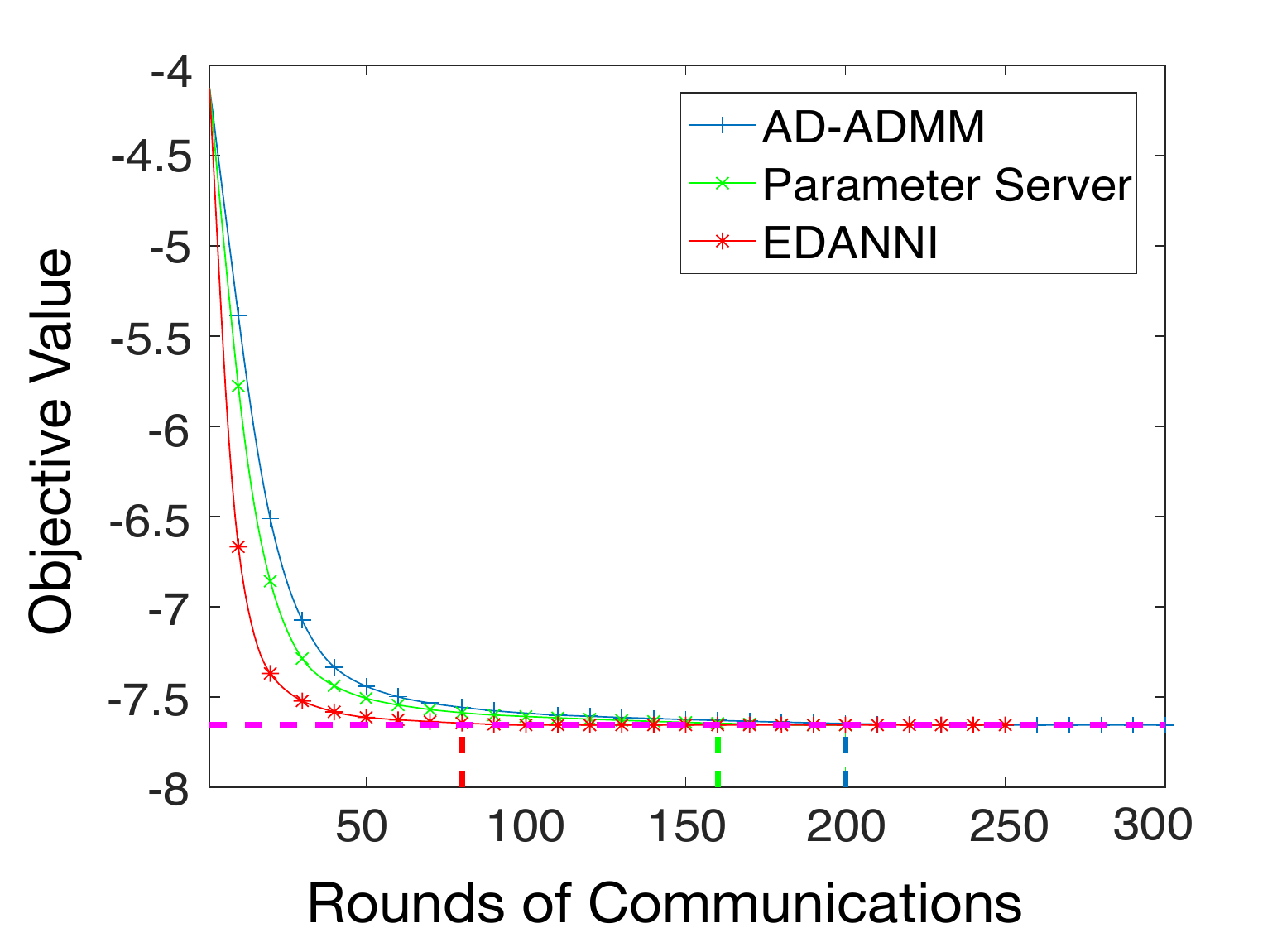}
		} 
		\hspace{-0.07in}\subfloat[{\scriptsize Asynchronous, \textbf{(3.2, 2.0, 1)}}]{
			\includegraphics[width=0.49\linewidth]{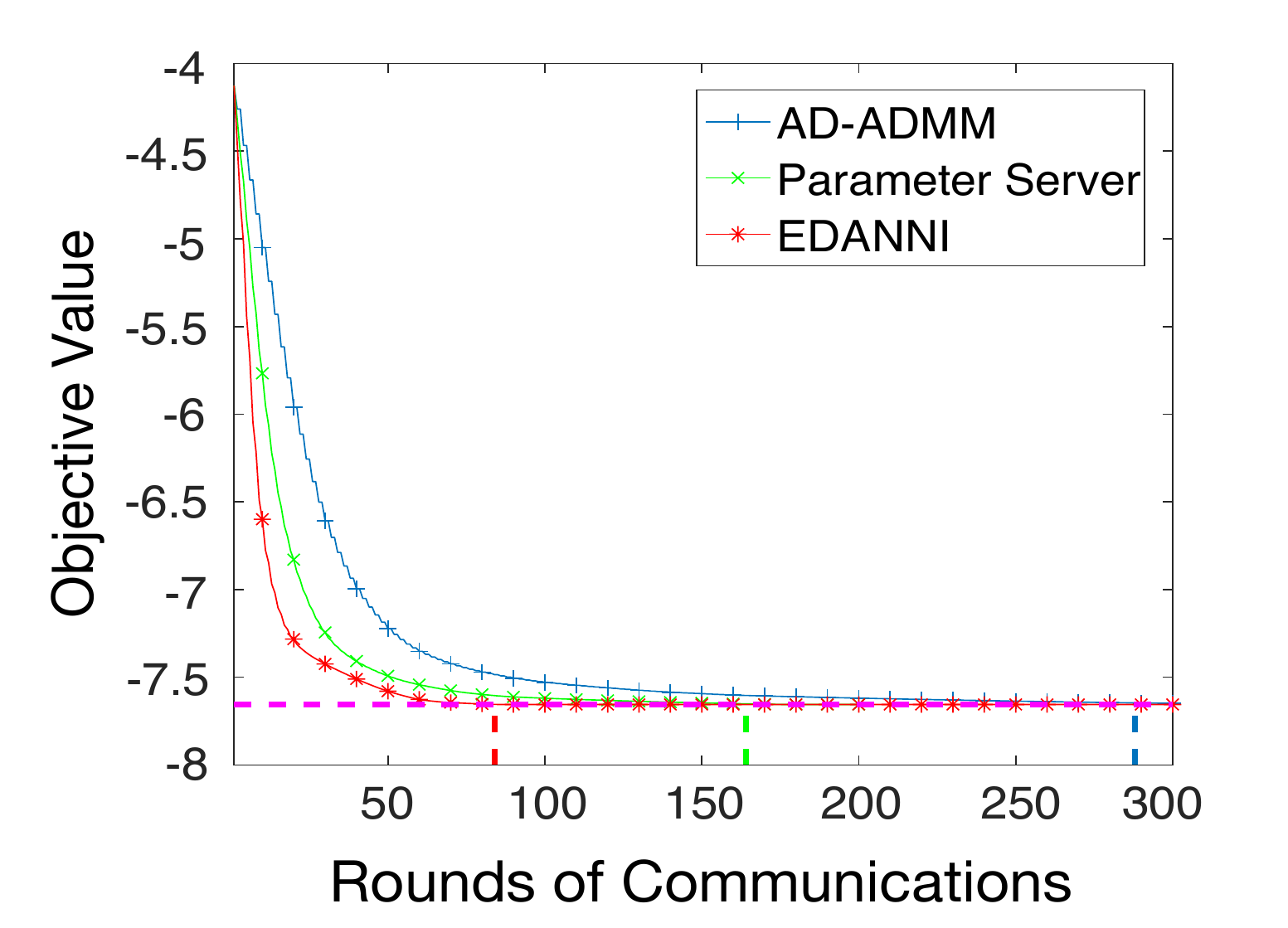}
		}
	\end{center}
	\vspace{-0.15in}	
	\caption{Comparison of candidate algorithms in sparse PCA. }\label{SPCAdiffalgr}
	\vspace{-0.12in}
\end{figure}

To verify the convergence conclusion of Theorem~\ref{mainthm} for nonconvex nonsmooth problems, we consider the following sparse PCA problem \cite{richtarik2012alternating}:
\begin{align}
\nonumber &\argmin_{\wb\in \RR^{p}} -\frac{1}{mn}\sum_{j\in [m]}\sum_{i \in [n]}\wb^TB_{ji}B_{ji}^T\wb + \theta \|\wb\|_1,\\
&\hspace{1.2em} \st \hspace{2em} \|\wb\| \leq 1
\end{align}
where $B_{ji} \in \RR^{p\times q}$ is a sparse matrix, $\forall j\in [m], ~~ i \in [n]$, and the regularization coefficient $\theta > 0$. Note that this is not a convex problem. In this example, we set  $m=3$,  $n=20$ $p=500$, $q = 1000$, and $\theta = 0.1$. Each matrix $B_{ji}\in \RR^{500\times1000}$ is a sparse random matrix with nearly $s=3000$ non-zero entries. The parameter $\rho$ in \eqref{update} is set to $\rho = 2\lambda_{\max}\left(\sum_{j\in[m]}B_{ji} B_{ji}^T\right)$. The candidate algorithms are compared in both the synchronous scenario and the asynchronous scenario. One can see from Figure~\ref{SPCAdiffalgr} (a) and Figure~\ref{SPCAdiffalgr} (b) that the proposed approach converges much faster than AD-ADMM and Parameter Server with much less communication cost. The results for other settings of $m, n, p, q, s$ in Table~\ref{table_spca} also verify such a conclusion.

\begin{figure}[!tbh]
	\centering 
	\vspace{0.01in}
	\begin{center}
		\hspace{-0.12in}\subfloat[{\scriptsize Convergence for different $\tau$'s}]{
			\includegraphics[width=0.475\linewidth]{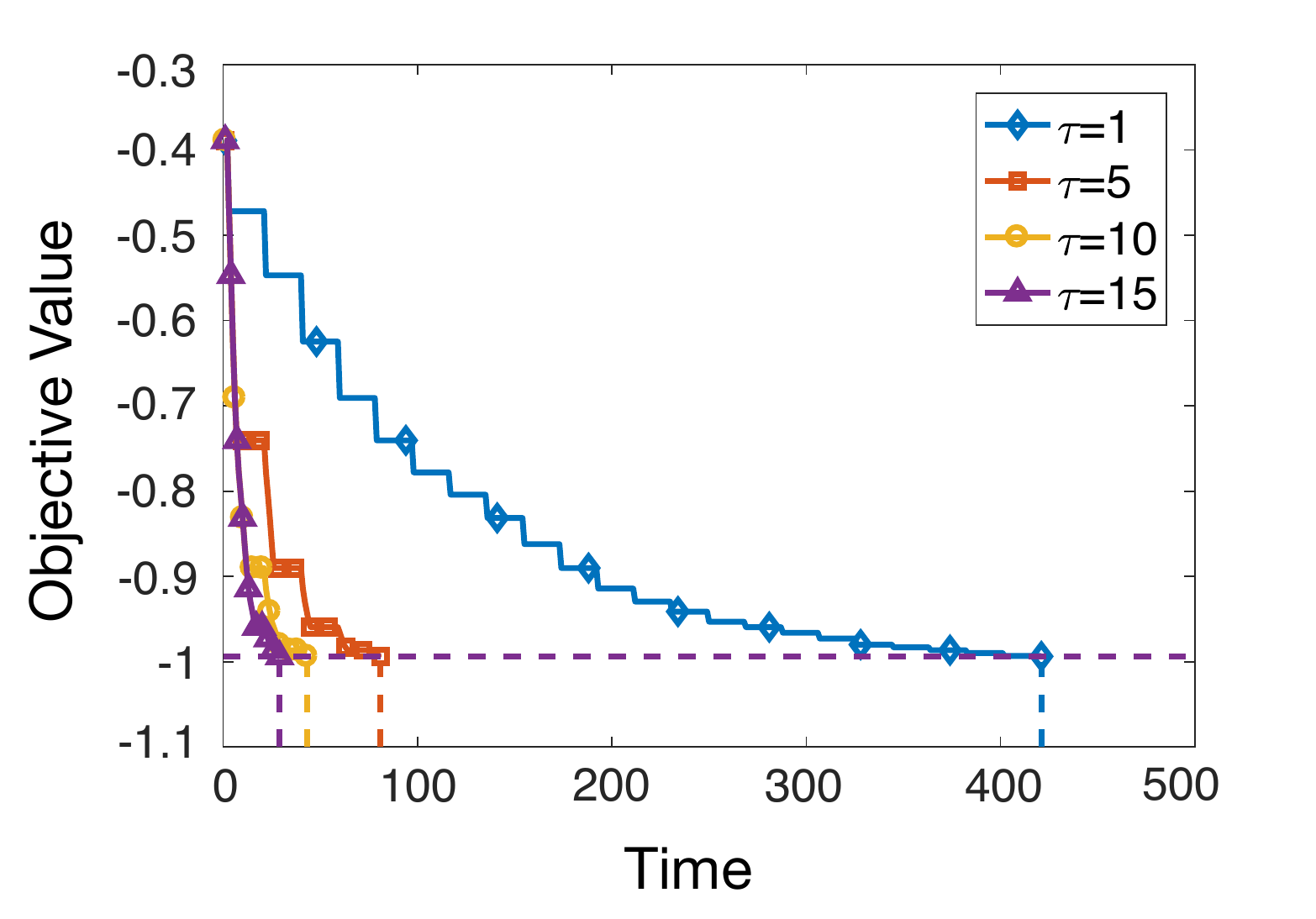}
		}
		\hspace{-0.06in}\subfloat[{\scriptsize Time table for different $\tau$'s}]{
			\includegraphics[width=0.485\linewidth]{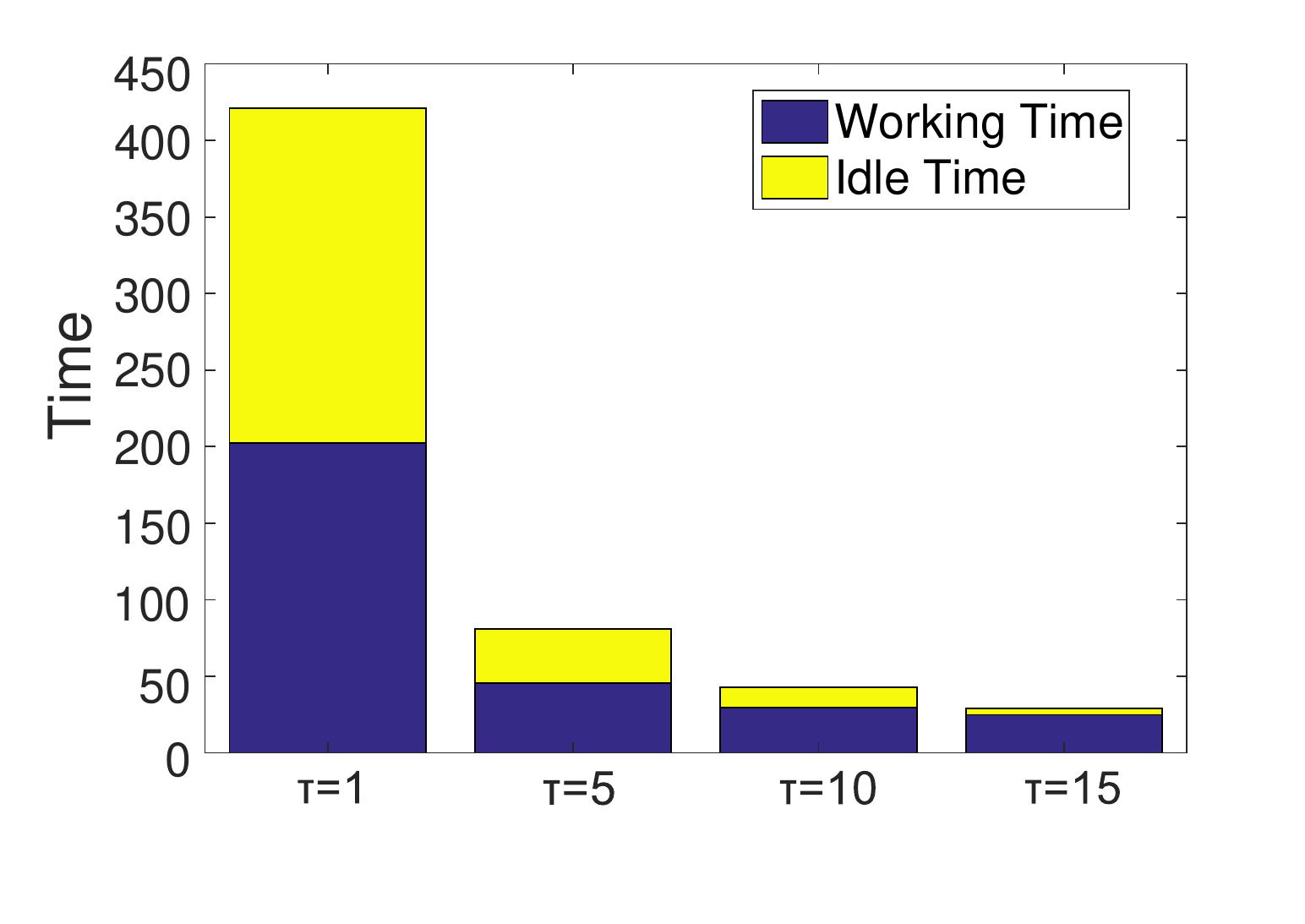}
		}
	\end{center}
	\vspace{-0.15in}	
	\caption{Comparison of EDANNI in sparse PCA with different settings of $\tau$ when $\rho=2$.   }\label{SPCAdifftau}
	\vspace{-0.2in}
\end{figure}

The performance of the proposed approach with different maximum tolerable delay $\tau$ is summarized in Figure~\ref{SPCAdifftau}. Here the distributed implementation is simulated in the same way as in the LASSO case. In Figure~\ref{SPCAdifftau} we set $m=6$,  $n=20$ $p=100$, $q = 1000$. One can observe from Figure~\ref{SPCAdifftau} (a) that in this example the convergence rate in terms of time is indeed affected by values of $\tau$. The running time when $\tau=15$ is much less than that when $\tau$ is small.  Similar to the LASSO case, Figure~\ref{SPCAdifftau} (b) shows that the ratio of the computing time in the overall running time increases closely to $1$ when the delay bound $\tau$ becomes $15$, therefore speeding up the convergence. Such results can be observed generally regardless of the choice of parameters $m, n, p, q, s, \theta$.

 \vspace{-0.2in}

\begin{table}[!tbh] \renewcommand*{\arraystretch}{1}
	\centering
	\caption{\small Comparison of the communication cost for SPCA}
	{\scriptsize \begin{tabular}{|p{1cm}|c|c|c|c|}
			\hline
			\multirow{2}{*}{\textbf{Type}} & \multirow{2}{*}{{\vspace{0.5em}\textbf{Parameters}}} & \multirow{2}{*}{\textbf{Method}}&\multirow{2}{*}{\textbf{Communication} }\\
			&(m,n,p,q,s,$\theta$)&& \\\hline
			\multirow{9}{*}{\hspace{-0.7em} \textbf{{\tiny Synchronous}}} &	\multirow{3}{*}{(20,100,50,100,50,0.1)}
			&\textbf{AD-ADMM}					&11.5\\\cline{3-4}
			&&\textbf{PS} &2.4\\\cline{3-4}
			&&\textbf{EDANNI}	   &\textbf{1}\\\cline{2-4}
			&	\multirow{3}{*}{(30,200,50,100,50,0.1)}
			&\textbf{AD-ADMM}				    	&30.0\\\cline{3-4}
			&&\textbf{PS} &2.5\\\cline{3-4}
			&&\textbf{EDANNI}	  &\textbf{1}\\\cline{2-4}
			&	\multirow{3}{*}{(3,20,500,1000,500,0.1)}
			&\textbf{AD-ADMM}					    &5.2\\\cline{3-4}
			&&\textbf{PS} &1.7\\\cline{3-4}
			&&\textbf{EDANNI}	   &\textbf{1}\\\hline
			\multirow{9}{*}{\hspace{-0.7em} \textbf{{\tiny Asynchronous}}} &	\multirow{2}{*}{(20,100,50,100,50,0.1)}
			&\textbf{AD-ADMM}					&12.5\\\cline{3-4}
			&\multirow{1}{*}{\hspace{0.1em} \textbf{{\tiny ($\tau = 3$)}}}&\textbf{PS} &4.1\\\cline{3-4}
			&&\textbf{EDANNI}	   &\textbf{1}\\\cline{2-4}
			&	\multirow{2}{*}{(30,200,50,100,50,0.1)}
			&\textbf{AD-ADMM}				    	&31.1\\\cline{3-4}
			&\multirow{1}{*}{\hspace{0.1em} \textbf{{\tiny ($\tau = 10$)}}}&\textbf{PS} &7.2\\\cline{3-4}
			&&\textbf{EDANNI}	   &\textbf{1}\\\cline{2-4}
			&	\multirow{2}{*}{(3,20,500,1000,500,0.1)}
			&\textbf{AD-ADMM}					    &6.3\\\cline{3-4}
			&\multirow{1}{*}{\hspace{0.1em} \textbf{{\tiny ($\tau = 3$)}}}&\textbf{PS} &2.0\\\cline{3-4}
			&&\textbf{EDANNI}	   &\textbf{1}\\\hline
	\end{tabular}}\label{table_spca}
\end{table}

\vspace{-0.1in}

%
%
%


\vspace{0.01in}
\section{Conclusion}
This paper proposes a communication-efficient distributed algorithm (EDANNI) solving a general problem \eqref{origin} in signal processing and machine learning under an asynchronous protocol. Theoretically, we prove the proposed algorithm converges to a stationary point in a sublinear rate, even in nonconvex nonsmooth scenarios. Moreover, unlike the previous work, linear convergence rate is established in strongly convex scenarios without any statistical assumptions of the local data. In experiments, we compare EDANNI with other state-of-the-art distributed algorithms in different applications, and the results show the superior performance of the proposed algorithm in terms of communication efficiency and the speed up caused by the asynchrony.

\newpage
\bibliographystyle{IEEEbib}
\bibliography{Renbib_jdh1}
\newpage
\appendices

\section{Proofs of Lemmata}
\vspace{0em}

\subsection{Proof of Lemma \ref{partdescent}}
Using optimality of $\xb^{t+1}$ in the update \eqref{update}, we have
\begin{align}
&{ \nonumber -\Bigg[\nabla\Lb_1(\xb^{t+1}) + \frac{1}{m} \sum_{j=1}^m \nabla\Lb_j(\xb^{t_j}) - \nabla\Lb_1(\xb^{t_1})\Bigg] }\\
&\hspace{1em} {\textstyle \in \partial\left[h(\xb) + \frac{\rho}{2}\left\|\xb - \xb^t\right\|^2 \right]. }
\end{align}
Recall that in Assumption~\ref{convsubproblem} (I), we define the convex modulus of $\frac{\rho}{2}\|\xb - \xb^{t}\|^2 + h(\xb) $ by $\gamma(\rho)$. It follows that

$\frac{\rho}{2}\|\xb^{t+1} - \xb^{t}\|^2 + h(\xb^{t+1}) - \left( \frac{\rho}{2}\|\xb^{t} - \xb^{t}\|^2 + h(\xb^{t})\right)$
\begin{align*}
&\leq \Bigg \langle { \nonumber -\Big[\nabla\Lb_1(\xb^{t+1}) + \frac{1}{m} \sum_{j=1}^m \nabla\Lb_j(\xb^{t_j}) - \nabla\Lb_1(\xb^{t_1})\Big] }, \\
& \hspace{1.2em}  \xb^{t+1} - \xb^{t} \Bigg \rangle - \frac{\gamma(\rho)}{2}\|\xb^{t+1} - \xb^t\|^2\\
&= -\Bigg \langle \nabla\Lb_1(\xb^{t+1}) + \frac{1}{m}\sum_{j=1}^m \nabla\Lb_j(\xb^{t_j})- \nabla\Lb_1(\xb^{t_1}),\\
& \Delta^{(t)} \Bigg \rangle  - \frac{\gamma(\rho)}{2}\big\|\Delta^{(t)}\big\|^2, 
\end{align*}
where we define $\Delta^{(t)}:= \xb^{t+1} - \xb^t$. 
Therefore
\begin{align}\label{insertinq}
\nonumber &\frac{\rho}{2}\|\xb^{t+1}-\xb^{t}\|^2 + h(\xb^{t+1}) - \frac{\rho}{2}\|\xb^{t}-\xb^{t-1}\|^2 - h(\xb^{t}) \\
\nonumber& = \frac{\rho}{2}\|\xb^{t+1}-\xb^{t}\|^2 + h(\xb^{t+1}) -  \frac{\rho}{2}\|\xb^{t}-\xb^{t}\|^2 - h(\xb^{t})  \\
\nonumber&+ \frac{\rho}{2}\|\xb^{t}-\xb^{t}\|^2 + h(\xb^{t})  -  \frac{\rho}{2}\|\xb^{t}-\xb^{t-1}\|^2 - h(\xb^{t}) \\
\nonumber&\leq -\Bigg \langle \nabla\Lb_1(\xb^{t+1}) + \frac{1}{m}\sum_{j=1}^m \nabla\Lb_j(\xb^{t_j})- \nabla\Lb_1(\xb^{t_1}),  \\
&  \Delta^{(t)} \Bigg \rangle - \frac{\gamma(\rho)}{2}\big\|\Delta^{(t)}\big\|^2 - \frac{\rho}{2}\big\|\Delta^{(t-1)}\big\|^2,
\end{align}
proving Lemma \ref{partdescent}. 

\subsection{Proof of Lemma \ref{descentlemma}}

It follows from Assumption~\ref{Lipdiff} that $\nabla\Lb_j(\xb)$ is Lipschitz continuous with constant $L$. Therefore we have 

\vspace{-0.3in}
\begin{align}\label{lipinq}
\nonumber&\frac{1}{m}\sum\limits_{j \in [m]} \Lb_j(\xb^{t+1}) - \frac{1}{m}\sum\limits_{j \in [m]} \Lb_j(\xb^{t})\\
\nonumber &\leq \left\langle \xb^{t+1} - \xb^t, \frac{1}{m}\sum_{j \in [m]}\nabla \Lb_j(\xb^{t})\right\rangle + \frac{L}{2}\|\xb^{t+1} - \xb^t\|^2\\
& = \left\langle \Delta^{(t)}, \frac{1}{m} \sum_{j\in [m]} \nabla \Lb_j(\xb^{t})\right\rangle + \frac{L}{2}\big\|\Delta^{(t)}\big\|^2.
\end{align}
%
By the above definition of fucntion $\Ft$, combining \eqref{insertinq} and \eqref{lipinq} results in
\begin{align}\label{Fdescent}
\nonumber &\text{F}(\xb^{t+1}, \xb^{t}) - \text{F}(\xb^{t}, \xb^{t-1}) \\
\nonumber &\overset{(b)}{\leq} -\Bigg\langle \nabla\Lb_1(\xb^{t+1}) + \frac{1}{m}\sum_{j \in [m]}\nabla \Lb_j(\xb^{t_j})-\nabla\Lb_1(\xb^{t_1}),\\ 
\nonumber & \hspace{1.1em} \Delta^{(t)} \Bigg\rangle -\frac{\gamma(\rho)}{2}\big\|\Delta^{(t)}\big\|^2 - \frac{\rho}{2}\big\|\Delta^{(t-1)}\big\|^2 \\
\nonumber &  + \left\langle\Delta^{(t)}, \frac{1}{m}\sum_{j \in [m]}\nabla \Lb_j(\xb^{t}) \right\rangle + \frac{L}{2}\big\|\Delta^{(t)}\big\|^2 \\
\nonumber & = -\Bigg\langle \nabla\Lb_1(\xb^{t+1}) + \frac{1}{m}\sum_{j \in [m]}\nabla \Lb_j(\xb^{t})-\nabla\Lb_1(\xb^{t}),  \\
\nonumber &  \hspace{1.1em} \Delta^{(t)} \Bigg\rangle - \frac{\gamma(\rho)}{2}\big\|\Delta^{(t)}\big\|^2 - \frac{\rho}{2}\big\|\Delta^{(t-1)}\big\|^2\\
\nonumber & +\left\langle\Delta^{(t)}, \frac{1}{m}\sum_{j \in [m]}\nabla \Lb_j(\xb^{t}) \right\rangle + \frac{L}{2}\big\|\Delta^{(t)}\big\|^2\\
\nonumber &  + \Bigg\langle\frac{1}{m}\sum_{j \in [m]}\nabla\Lb_j(\xb^{t})- \frac{1}{m}\sum_{j \in [m]}\nabla \Lb_j(\xb^{t_j}),\Delta^{(t)}\Bigg\rangle\\
\nonumber & +\left \langle\nabla\Lb_1(\xb^{t_1}) - \nabla\Lb_1(\xb^{t}), \Delta^{(t)} \right \rangle\\
\nonumber & \leq \left(\frac{3L}{2} - \frac{\gamma(\rho)}{2}\right)\big\|\Delta^{(t)}\big\|^2 - \frac{\rho}{2}\big\|\Delta^{(t-1)}\big\|^2\\
\nonumber & +\underbrace{\left\langle\frac{1}{m}\sum_{j \in [m]}\nabla\Lb_j(\xb^{t})- \frac{1}{m}\sum_{j \in [m]}\nabla \Lb_j(\xb^{t_j}),\Delta^{(t)}\right\rangle }_{(P1)}\\
& \underbrace{ +\left \langle\nabla\Lb_1(\xb^{t_1}) - \nabla\Lb_1(\xb^{t}), \Delta^{(t)} \right \rangle}_{(P1)}, 
\end{align}
where inequality $(b)$ is due to Lemma~\ref{partdescent} and Assumption~\ref{Lipdiff}. 

Note that

$\nabla\Lb_1(\xb^{t_1}) - \nabla\Lb_1(\xb^{t})$
\begin{align*}
\hspace{3em} = \sum_{k = 1}^{t-t_1}\big(\nabla\Lb_1(\xb^{t-k})-\nabla\Lb_1(\xb^{t-k+1})\big),
\end{align*}
which implies

$\left\|\nabla\Lb_1(\xb^{t_1}) - \nabla\Lb_1(\xb^{t})\right\|$
\begin{align*}
&\leq \sum_{k = 1}^{t-t_1}\left\|\nabla\Lb_1(\xb^{t-k})-\nabla\Lb_1(\xb^{t-k+1})\right\|\\
& \leq \sum_{k = 1}^{t-t_1} L \left\|\xb^{t-k}- \xb^{t-k+1}\right\|\\
& \leq \sum_{k = 1}^{\tau} L \left\|\xb^{t-k}- \xb^{t-k+1}\right\|\\
& = \sum_{k = 1}^{\tau} L \big\|\Delta^{(t-k)}\big\|. 
\end{align*}
Similarly, we can see

$\left\| \frac{1}{m}\sum_{j \in [m]} \nabla \Lb_j(\xb^{t}) - \frac{1}{m}\sum_{j \in [m]} \nabla \Lb_j(\xb^{t_j}) \right\|$
\begin{align*}
\hspace{8em} \leq  \sum_{k = 1}^{\tau} L \big\|\Delta^{(t-k)}\big\|.
\end{align*}
These two inequalities result in
\begin{align}\label{2bound}
\nonumber(P1) &\leq 2L\sum_{k = 1}^{\tau} \big\|\Delta^{(t-k)}\big\|\big\|\Delta^{(t)}\big\|\\
\nonumber&\leq L\sum_{k = 1}^\tau \left(\frac{1}{\delta}\big\|\Delta^{(t-k)}\big\|^2 + \delta\big\| \Delta^{(t)}\big\|^2\right)\\
&\leq \frac{L}{\delta}\sum_{k = 1}^\tau \big\|\Delta^{(t-k)}\big\|^2 + L\delta\tau\big\| \Delta^{(t)}\big\|^2,
\end{align}
where in the second inequality we apply the fact that
\begin{align*}
a\cdot b\leq \frac{1}{2}\left(\frac{1}{\delta}a^2 + \delta b^2\right)
\end{align*}
for any $a, b, \delta > 0$. By inserting \eqref{2bound} into \eqref{Fdescent} we have

$\text{F}(\xb^{t+1}, \xb^{t}) - \text{F}(\xb^{t}, \xb^{t-1})$
\begin{align*}
\nonumber\leq \left(\frac{3L}{2} - \frac{\gamma(\rho)}{2}+L\delta\tau\right)\big\|\Delta^{(t)}\big\|^2 - \frac{\rho}{2}\big\|\Delta^{(t-1)}\big\|^2 \\
+ \frac{L}{\delta}\sum_{k = 1}^{\tau} \big\|\Delta^{(t-k)}\big\|^2,
\end{align*}
proving the conclusion of Lemma \ref{descentlemma}. 

\subsection{Proof of Lemma \ref{Flowbound}}
First of all, summing the above inequality \eqref{adjacdescent1} of Lemma~\ref{descentlemma} over $t$ yields

$\text{F}(\xb^{T+1}, \xb^{T}) - \text{F}(\xb^{1}, \xb^{0})$
\begin{align*}
\leq \sum_{t=0}^T \left(\frac{3L}{2} - \frac{\gamma(\rho)}{2} +L\delta\tau\right)\big\| \Delta^{(t)}\big\|^2 \\
+ \sum_{t=0}^T\left( \frac{L\tau}{\delta} -\frac{\rho}{2}\right)\big\| \Delta^{(t-1)}\big\|^2. 
\end{align*}
If $\rho$ satisfies Assumption \ref{convsubproblem}, then it holds that
\begin{align*}
\text{F}(\xb^{T+1}, \xb^{T}) - \text{F}(\xb^{1}, \xb^{0}) < 0.
\end{align*}	
By taking $\xb^{T}$ as the initial point, similarly we have
\begin{align*}
\text{F}(\xb^{2T+1}, \xb^{2T}) - \text{F}(\xb^{T+1}, \xb^{T}) < 0.
\end{align*}
Continuing this process we get a decreasing subsequence $\left\{\text{F}(\xb^{kT+1}, \xb^{kT})\right\}_{k=0,1,\cdots}$. Therefore there exists a constant $\bar{F}_0$ such that
\begin{align}\label{Flimit} 
\Ft(\xb^{kT+1}, \xb^{kT}) \leq \bar{F}_0. 
\end{align}
When starting with $\xb^1, \cdots, \xb^{T-1}$, with similar analysis we can prove that there exists constants $\bar{F}_1, \cdots, \bar{F}_{T-1}$ such that 
\begin{align}\label{Flimit} 
\Ft(\xb^{kT+l+1}, \xb^{kT+l}) \leq \bar{F_l}, \hspace{1.1em} 
\end{align}
for $l = 1,2,\cdots,T-1$. 
Define $\bar{\Ft} := \max\left\{ \bar{F}_0, \cdots, \bar{F}_{T-1}\right\}$, then 
\begin{align*}
\text{F}(\xb^{t+1}, \xb^{t}) < \bar{\Ft} < +\infty, \hspace{1.1em} \forall t \in \NN. 
\end{align*}

On the other hand, let $\underline{\Ft} := \underline{L}$, then by the definition of $\Ft(\xb^{t+1}, \xb^{t})$ and Assumption~\ref{convsubproblem}, we have

$\text{F}(\xb^{t+1}, \xb^{t}) $
\begin{align*}
&= \frac{1}{m}\sum_{j\in [m]}\Lb_j(\xb^{t+1}) + \frac{\rho}{2}\left\| \xb^{t+1}- \xb^t\right\|^2 + h(\xb^{t+1})\\
&\geq \frac{1}{m}\sum_{j\in [m]}\Lb_j(\xb^{t+1}) + h(\xb^{t+1}) \\
&=  \Lb(\xb^{t+1}) >\underline{L} = \underline{\Ft} > -\infty,
\end{align*}
for any $t\in \NN$.
Therefore the boundedness of function $\Ft$ in Lemma \ref{Flowbound} is proved. 
%
\vspace{1em}

\subsection{Proof of Lemma~\ref{upbound_optm}}
To prove the convergence rate, we first need to bound $\left(\Ft(\xb^{t+1}, \xb^t) - \Ft(\xb^*,\xb^*) \right)$, where $\Ft(\xb, \xb^t) = \frac{1}{m} \sum\limits_{j\in[m]} \Lb_j(\xb) + \frac{\rho}{2}\|\xb -\xb^t\|^2 + h(\xb) \geq \Ft(\xb^*, \xb^*) = \frac{1}{m} \sum\limits_{j\in[m]} \Lb_j(\xb^*) + h(\xb^*)$.

By the optimality of $\xb^{t+1}$ in the update \eqref{update}, we have
\begin{align}\label{firstorder_optimality}
\nonumber	&\Big( \nabla \Lb_1(\xb^{t+1}) + \frac{1}{m}\sum_{j\in[m]}\nabla \Lb_j(\xb^{t_j})-\nabla \Lb_1(\xb^{t})\\
&\textstyle + \partial h(\xb^{t+1}) + \rho(\xb^{t+1}-\xb^{t})\Big)^\top (\xb^{t+1}-\xb) \leq 0,
\end{align}
for all $\xb\in \RR^p$. Letting $\xb = \xb^*$ implies
\begin{align}\label{firstorder_optimality_xstar}
\nonumber & \Big( \nabla \Lb_1(\xb^{t+1}) + \frac{1}{m}\sum_{j\in[m]}\nabla \Lb_j(\xb^{t_j})-\nabla \Lb_1(\xb^{t}) \\
&\textstyle + \partial h(\xb^{t+1}) + \rho(\xb^{t+1}-\xb^{t})\Big)^\top (\xb^{t+1}-\xb^*) \leq 0. 
\end{align}

By the strong convexity of $\Lb_j$ one has
\begin{align}\label{strong_convexity} 
\nonumber\Lb_j(\yb) \geq \Lb_j(\xb) + \left(\nabla \Lb_j(\xb)\right)^\top(\yb - \xb) + \frac{\sigma^2}{2} \|\yb-\xb\|^2,\\
\forall \xb, ~ \yb \in \RR^p.
\end{align}

Setting $\yb=\xb^*$, $\xb=\xb^{t+1}$ in \eqref{strong_convexity} we have
\begin{align*}
\nonumber \Lb_j(\xb^*) &\geq \Lb_j(\xb^{t+1}) + \left(\nabla \Lb_j(\xb^{t+1})\right)^\top(\xb^* - \xb^{t+1}) \\
&+ \frac{\sigma^2}{2} \|\xb^*-\xb^{t+1}\|^2,
\end{align*}
which further implies that 
\begin{align}\label{strong_convexity_x}
\big(\nabla \nonumber &\Lb_j(\xb^{t+1})\big)^\top(\xb^{t+1}-\xb^*) \\
&\geq \Lb_j(\xb^{t+1}) -\Lb_j(\xb^*) + \frac{\sigma^2}{2} \|\xb^*-\xb^{t+1}\|^2.
\end{align}

Note that \eqref{firstorder_optimality_xstar} implies that
\begin{align}\label{optimal_result}
&\nonumber \Big(\nabla \Lb_1(\xb^{t+1}) - \nabla \Lb_1(\xb^{t}) + \Big(\frac{1}{m} \sum_{j\in [m]}\nabla \Lb_j(\xb^{t_j}) \\
\nonumber &- \frac{1}{m} \sum_{j\in [m]} \nabla \Lb_j(\xb^{t+1})\Big) +\frac{1}{m}\sum_j \nabla\Lb_j(\xb^{t+1})\\
&+\partial h(\xb^{t+1}) + \rho(\xb^{t+1} - \xb^t)\Big)^\top (\xb^{t+1}-\xb^*) \leq 0. 
\end{align}
Summing \eqref{strong_convexity_x} over $j\in [m]$ gives that
\begin{align}\label{strong_convexity_result}
\nonumber\frac{1}{m}\sum_{j\in[m]}\left(\nabla \Lb_j(\xb^{t+1})\right)^\top(\xb^{t+1}-\xb^*) \geq \frac{1}{m}\sum_{j\in[m]}\Lb_j(\xb^{t+1})\\ -\frac{1}{m}\sum_{j\in[m]}\Lb_j(\xb^*) + \frac{\sigma^2}{2} \|\xb^*-\xb^{t+1}\|^2.
\end{align}
Putting \eqref{strong_convexity_result} into \eqref{optimal_result}, we have
\begin{align}\label{combine1}
&\nonumber \Big( \nabla \Lb_1(\xb^{t+1}) - \nabla \Lb_1(\xb^{t}) + \Big(\frac{1}{m} \sum_{j\in [m]}\nabla \Lb_j(\xb^{t_j}) \\
&\nonumber - \frac{1}{m} \sum_{j\in [m]} \nabla \Lb_j(\xb^{t+1})\Big)\Big)^\top(\xb^{t+1}- \xb^*) \\ 
&\nonumber + \Big(\frac{1}{m}\sum_{j\in [m]}\Lb_j(\xb^{t+1}) - \frac{1}{m}\sum_{j\in [m]}\Lb_j(\xb^*)\Big) \\
\nonumber &+\frac{\sigma^2}{2}\|\xb^* - \xb^{t+1}\|^2 + \partial h(\xb^{t+1})(\xb^{t+1} - \xb^*) \\
&+ \rho(\xb^{t+1} - \xb^t)^\top (\xb^{t+1} - \xb^*) \leq 0. 
\end{align}

Since $h(\xb)$ is convex, we have
\begin{align*}
h(\xb^{t+1}) - h(\xb^*) \leq \partial h(\xb^{t+1})(\xb^{t+1} - \xb^*).
\end{align*}
Putting it into \eqref{combine1}, we have
\begin{align}\label{combine2}
\nonumber &\Big( \nabla \Lb_1(\xb^{t+1}) - \nabla \Lb_1(\xb^{t}) + \Big(\frac{1}{m} \sum_{j\in [m]}\nabla \Lb_j(\xb^{t_j})\\
\nonumber &- \frac{1}{m} \sum_{j\in [m]} \nabla \Lb_j(\xb^{t+1})\Big)\Big)^\top (\xb^{t+1}-\xb^*) \\
\nonumber &+ \Big(\frac{1}{m} \sum_{j\in [m]}\Lb_j(\xb^{t+1}) - \frac{1}{m} \sum_{j\in [m]} \Lb_j(\xb^*)\Big)\\
\nonumber &+\frac{\sigma^2}{2}\|\xb^*-\xb^{t+1}\|^2 + h(\xb^{t+1}) - h(\xb^*) \\
&+ \rho(\xb^{t+1} - \xb^t)^\top (\xb^{t+1} - \xb^*) \leq 0.
\end{align}\vspace{-0.3in}

Note that
\begin{align}\label{polar}
\nonumber\rho (\xb^{t+1} - \xb^{t})^\top(\xb^{t+1} - \xb^*) = \frac{\rho}{2}\|\xb^{t+1}-\xb^*\|^2 \\
+\frac{\rho}{2}\|\xb^{t+1}-\xb^{t}\|^2- \frac{\rho}{2}\|\xb^{t}-\xb^*\|^2.
\end{align}\vspace{-0.3in}

\hspace{-0.75em}Then putting \eqref{polar} into \eqref{combine2}, one obtains
\begin{align*}
\Big( \nabla \Lb_1(\xb^{t+1}) - \nabla \Lb_1(\xb^{t}) + \Big(\frac{1}{m} \sum_{j\in [m]}\nabla \Lb_j(\xb^{t_j}) \hspace{3.3em} \\
- \frac{1}{m} \sum_{j\in [m]} \nabla \Lb_j(\xb^{t+1})\Big)\Big)^\top(\xb^{t+1}-\xb^*)\hspace{3.3em} \\
+ \Big(\frac{1}{m} \sum_{j\in [m]}\Lb_j(\xb^{t+1}) - \frac{1}{m} \sum_{j\in [m]} \Lb_j(\xb^*)\Big)\hspace{3.3em}\\
+ \frac{\sigma^2}{2}\|\xb^*-\xb^{t+1}\|^2 + h(\xb^{t+1}) - h(\xb^*)\hspace{3.3em} \\
+\frac{\rho}{2}\|\xb^{t+1}-\xb^t\|^2+ \frac{\rho}{2}\|\xb^{t+1}-\xb^*\|^2 - \frac{\rho}{2}\|\xb^{t}-\xb^*\|^2\leq 0,\hspace{1.6em}  
\end{align*}
which is equivalent to
\begin{align}\label{eqn101}
\nonumber&\Ft(\xb^{t+1}, \xb^t) - \Ft(\xb^*, \xb^*) \leq -\frac{\sigma^2}{2}\|\xb^*-\xb^{t+1}\|^2\\
\nonumber &- \Big( \nabla \Lb_1(\xb^{t+1}) - \nabla \Lb_1(\xb^{t}) + \Big(\frac{1}{m} \sum_{j\in [m]}\nabla \Lb_j(\xb^{t_j}) \\
\nonumber & - \frac{1}{m} \sum_{j\in [m]} \nabla \Lb_j(\xb^{t+1})\Big)\Big)^\top (\xb^{t+1} - \xb^*) \\
& -\frac{\rho}{2}\|\xb^{t+1}-\xb^*\|^2 + \frac{\rho}{2}\|\xb^t - \xb^*\|^2.
\end{align}
Now, note that
\begin{align}\label{gradientinsert}
\nonumber &\Big( \nabla \Lb_1(\xb^{t+1}) - \nabla \Lb_1(\xb^{t}) + \Big(\frac{1}{m} \sum_{j\in [m]}\nabla \Lb_j(\xb^{t_j}) \\
\nonumber& - \frac{1}{m} \sum_{j\in [m]} \nabla \Lb_j(\xb^{t+1})\Big)\Big)^\top(\xb^{t+1}-\xb^*)\\
\nonumber &=\Big( \nabla \Lb_1(\xb^{t+1}) - \nabla \Lb_1(\xb^{t}) + \frac{1}{m} \sum_{j\in [m]}\nabla \Lb_j(\xb^{t}) \\
\nonumber& - \frac{1}{m} \sum_{j\in [m]} \nabla \Lb_j(\xb^{t+1})\Big)^\top(\xb^{t+1}-\xb^*)\\
&\underbrace{+ \Big( \frac{1}{m} \sum_{j\in [m]}\nabla \Lb_j(\xb^{t_j})-\frac{1}{m} \sum_{j\in [m]}\nabla \Lb_j(\xb^{t})\Big)^\top(\xb^{t+1}-\xb^*)}_{(P5)}.
\end{align}

By the Mean Value Theorem one has 
\begin{align*}
\nonumber &\Big( \nabla \Lb_1(\xb^{t+1}) - \nabla \Lb_1(\xb^{t}) + \Big(\frac{1}{m} \sum_{j\in [m]}\nabla \Lb_j(\xb^{t}) \\
\nonumber& - \frac{1}{m} \sum_{j\in [m]} \nabla \Lb_j(\xb^{t+1})\Big)\Big)^\top(\xb^{t+1}-\xb^*)\\
\nonumber &= (\xb^{t+1} - \xb^t)^\top \Big[\nabla^2\Lb_1(\xi) - \frac{1}{m}\sum_{j\in [m]} \nabla^2\Lb_j(\xi) \Big]\\
\nonumber& \hspace{1em}\cdot(\xb^{t+1} - \xb^*)\\
\nonumber&= \frac{1}{2}\|\xb^{t+1}-\xb^*\|_\Sigma^2 - \frac{1}{2}\|\xb^{t+1}-\xb^{t+1}\|_\Sigma^2 \\
&+\frac{1}{2}\|\xb^{t}-\xb^{t+1}\|_\Sigma^2 - \frac{1}{2}\|\xb^{t}-\xb^*\|_\Sigma^2,
\end{align*}
where $\Sigma:= \nabla^2\Lb_1(\xi) - \frac{1}{m}\sum_{j\in [m]} \nabla^2\Lb_j(\xi)$. 
It follows that 
\begin{align}
\nonumber &\Bigg( \nabla \Lb_1(\xb^{t+1}) - \nabla \Lb_1(\xb^{t}) + \Big(\frac{1}{m} \sum_{j\in [m]}\nabla \Lb_j(\xb^{t}) \\
\nonumber & -\frac{1}{m} \sum_{j\in [m]} \nabla \Lb_j(\xb^{t+1})\Big)\Bigg)^\top(\xb^{t+1}-\xb^*)\\
\nonumber&= \frac{1}{2}\|\xb^{t+1}-\xb^*\|_\Sigma^2 +\frac{1}{2}\|\xb^{t}-\xb^{t+1}\|_\Sigma^2 \\
\nonumber&- \frac{1}{2}\|\xb^{t}-\xb^{t+1}+\xb^{t+1}-\xb^*\|_\Sigma^2
\end{align}
\begin{align} \label{eqn16}
\nonumber&\geq \frac{1}{2}\|\xb^{t+1}-\xb^*\|_\Sigma^2 +\frac{1}{2}\|\xb^{t}-\xb^{t+1]}\|_\Sigma^2 \\
\nonumber &- \frac{1}{2}(1+\delta_1)\|\xb^{t+1}-\xb^t\|_\Sigma^2 - \frac{1}{2}(1+1/\delta_1)\|\xb^{t+1}-\xb^*\|_\Sigma^2\\
&\geq -\frac{1}{2}\delta_1\|\xb^{t}-\xb^{t+1}\|_\Sigma^2 - \frac{1}{2\delta_1}\|\xb^{t+1}-\xb^*\|_\Sigma^2.
\end{align}
Note that 
\begin{align*}\label{eqn17}
& \frac{\rho}{2}\|\xb^t - \xb^*\|^2 = \frac{\rho}{2}\|\xb^t -\xb^{t+1} + \xb^{t+1}- \xb^*\|^2\\
&\leq \frac{\rho}{2}(1+\delta_1)\|\xb^t -\xb^{t+1}\|^2 + \frac{\rho}{2}(1+1/\delta_1) \|\xb^{t+1}- \xb^*\|^2, 
\end{align*}
which implies that
\begin{align}
\nonumber &-\frac{\rho}{2}\|\xb^{t+1}- \xb^*\|^2 + \frac{\rho}{2}\|\xb^t - \xb^*\|^2\\
& \leq \frac{\rho}{2}(1+\delta_1)\|\xb^t -\xb^{t+1}\|^2 + \frac{\rho}{2}\frac{1}{\delta_1} \|\xb^{t+1}- \xb^*\|^2. 
\end{align}\vspace{-0.6em}

Putting \eqref{gradientinsert}, \eqref{eqn16}, and \eqref{eqn17} into \eqref{eqn101}, we can bound the optimality gap of function $\Ft$ by\vspace{-0.4em}
\begin{align}\label{descentF_31}
&\nonumber \Ft(\xb^{t+1}, \xb^t) - \Ft(\xb^*, \xb^*) \leq -\frac{\sigma^2}{2}\|\xb^*-\xb^{t+1}\|^2 \\
&\nonumber +\frac{1}{2}\delta_1\|\xb^{t+1} - \xb^t\|_\Sigma^2 + \frac{1}{2\delta_1}\|\xb^{t+1} - \xb^*\|_\Sigma^2\\
&\nonumber -\frac{\rho}{2}\|\xb^{t+1}- \xb^*\|^2 + \frac{\rho}{2}\|\xb^t - \xb^*\|^2-(P5)\\
%
&\nonumber = -\frac{\sigma^2}{2}\|\xb^*-\xb^{t+1}\|^2  + \frac{1}{2\delta_1}\|\xb^{t+1} - \xb^*\|_\Sigma^2\\
& \nonumber  +\frac{\rho}{2\delta_1}\|\xb^{t+1} - \xb^*\|^2 +\frac{1}{2}\delta_1\|\xb^{t+1} - \xb^t\|_\Sigma^2 \\
& + \frac{\rho}{2}(1+\delta_1)\|\xb^t -\xb^{t+1}\|^2 - (P5). 
\end{align}

Now we bound (P5) on the RHS of \eqref{descentF_31}. For some $\delta_1>0$ one has
\begin{align*}
&(P5) \\ 
&\geq -\frac{\delta_1}{2m} \sum_{j\in[m]}\|\nabla\Lb_j(\xb^{t_j})- \nabla\Lb_j(\xb^{t})\|^2 - \frac{1}{2\delta_1}\|\xb^{t+1}- \xb^{*}\|^2\\
&\geq -\frac{\delta_1}{2m}L^2 \sum_{j\in[m]}\|\xb^{t_j}- \xb^{t}\|^2 - \frac{1}{2\delta_1}\|\xb^{t+1}- \xb^{*}\|^2.
\end{align*}
Therefore we have
\begin{align}\label{eqn10}
&\nonumber \Ft(\xb^{t+1}, \xb^t) - \Ft(\xb^*, \xb^*)\\
&\nonumber \leq -\frac{\sigma^2}{2}\|\xb^*-\xb^{t+1}\|^2  + \frac{1}{2\delta_1}\|\xb^{t+1} - \xb^*\|_\Sigma^2\\
& \nonumber  +\frac{\rho}{2\delta_1}\|\xb^{t+1} - \xb^*\|^2 +\frac{1}{2}\delta_1\|\xb^{t+1} - \xb^t\|_\Sigma^2 \\
&\nonumber + \frac{\rho}{2}(1+\delta_1)\|\xb^t -\xb^{t+1}\|^2 \\
& +\frac{\delta_1}{2m}L^2 \sum_{j\in[m]}\|\xb^{t_j}- \xb^{t}\|^2 + \frac{1}{2\delta_1}\|\xb^{t+1}- \xb^{*}\|^2. 
\end{align}

Let $\frac{\sigma^2}{2} > \frac{2L}{2\delta_1} + \frac{\rho}{2\delta_1} + \frac{1}{2\delta_1}$, i.e., $\sigma^2> \frac{2L}{\delta_1} + \frac{\rho}{\delta_1}+ \frac{1}{\delta_1}$, then
\begin{align}
\nonumber&\Ft(\xb^{t+1}, \xb^t) - \Ft(\xb^*, \xb^*) \\
\nonumber &\leq \frac{1}{2}\delta_1\|\xb^{t+1} - \xb^t\|_\Sigma^2 + \frac{\rho}{2}(1+\delta_1)\|\xb^t -\xb^{t+1}\|^2\\
\nonumber & +\frac{\delta_1}{2m}L^2 \sum_{j\in[m]}\|\xb^{t_j}- \xb^{t}\|^2\\
\nonumber &\leq \delta_1 L\|\xb^{t+1} - \xb^t\|^2 + \frac{\rho}{2}(1+\delta_1)\|\xb^t -\xb^{t+1}\|^2\\
& +\frac{\delta_1}{2m}L^2 \sum_{j\in[m]}\|\xb^{t_j}- \xb^{t}\|^2,
\end{align}
from which we have
\begin{align*}
\nonumber& \frac{1}{\frac{\rho}{2}(1+\delta_1)+\delta_1}\left( \Ft(\xb^{t+1}, \xb^t) - \Ft(\xb^*, \xb^*) \right)\\
\nonumber &\leq \frac{\delta_1 L}{\frac{\rho}{2}(1+\delta_1)+\delta_1}\|\xb^t-\xb^{t+1}\|^2 \\
\nonumber&+ \frac{\frac{\rho}{2}(1+\delta_1)}{\frac{\rho}{2}(1+\delta_1)+\delta_1}\|\xb^t - \xb^{t+1} \|^2 \\
& + \frac{1}{\frac{\rho}{2}(1+\delta_1)+\delta_1}\frac{\delta_1}{2m}L^2 \sum_{j\in[m]}\|\xb^{t_j}- \xb^{t}\|^2,
\end{align*}
therefore proving Lemma~\ref{upbound_optm}.

\end{document}